\documentclass[3p,times,sort]{elsarticle}

\journal{Applied Mathematics and Computation}

\usepackage{graphicx}
\usepackage[reqno]{amsmath}
\usepackage{amssymb}
\let\vec\relax 
\usepackage{epstopdf}
\usepackage{amsthm}
\let\vec\relax 
\usepackage[]{graphicx}               
\usepackage{bm}
\usepackage{esvect}
\usepackage{ulem}
\usepackage{adjustbox}
\usepackage{accents}
\usepackage{stmaryrd}
\usepackage{xcolor}
\usepackage{enumitem}
\usepackage{hyperref}
\usepackage{booktabs}
\usepackage{enumitem}
\usepackage{longtable}

\usepackage{xargs}    

\usepackage[colorinlistoftodos,prependcaption,textsize=tiny]{todonotes}
\newcommandx{\unsure}[2][1=]{\todo[linecolor=blue,backgroundcolor=blue!25,bordercolor=blue,#1]{#2}}
\newcommandx{\changeThis}[2][1=]{\todo[linecolor=red,backgroundcolor=red!25,bordercolor=red,#1]{#2}}

\newtheorem{theorem}{Theorem}[section]

\theoremstyle{definition}
\newtheorem{example}[theorem]{Example}

\theoremstyle{remark}
\newtheorem{remark}[theorem]{Remark}

\numberwithin{equation}{section}

\begin{document}


\definecolor{chcol}{rgb}{0.4,0.,0.9}
\newcommand{\change}[1]{\textcolor{chcol}{#1}}
\newcommand{\viscosity}{\mu}     
\renewcommand{\Pr}{\text{Pr}}     
\renewcommand{\Re}{\text{Re}}   
\newcommand{\Ma}{\text{Ma}}     
\newcommand{\halfComma}{\kern 0.083334em}
\newcommand{\pderivative}[2]{\frac{\partial #1}{\partial #2}}
\newcommand{\ddt}{\frac{d}{dt}}
\newcommand{\avg}[1]{\left\{\hspace*{-3pt}\left\{#1\right\}\hspace*{-3pt}\right\}}
\newcommand{\jump}[1]{\ensuremath{\left\llbracket #1 \right\rrbracket}}
\newcommand{\shat}{\ensuremath{\hat{s}}}        
\newcommand{\dS}{{\,\operatorname{dS}}}         
\newcommand{\PBT}{{\,\operatorname{PBT}}}    
\newcommand{\ec}{{\mathrm{EC}}}                     
\newcommand{\es}{{\mathrm{ES}}}                     
\newcommand{\ent}{{S}}                                     
\newcommand{\ma}{-}                                        
\renewcommand{\sl}{+}                                      
\newcommand{\dx}{\,\text{d}x}                           
\newcommand\iprod[1]{\left\langle #1\right\rangle}                                             
\newcommand\inorm[1]{\left |\left| #1\right|\right|}                                               
\newcommand\iprodN[1]{\left\langle #1\right\rangle_{\!N}}                                 
\newcommand\inormN[1]{\left |\left| #1\right|\right|_{N}}                                     
\newcommand\irefInt{\int\limits_{-1}^{1}}                                                            
\newcommand\ivolN{\int\limits_{ N }\! }                                                              
\newcommand\isurfE{\int\limits_{\partial E }\! }                                                   
\newcommand\volEN{\mkern-11mu\int\limits_{E , N}\mkern-5mu }                    
\newcommand\isurfEN{\mkern-11mu\int\limits_{\partial E , N}\mkern-11mu }    
\newcommand\vnorm[1]{\left| #1\right|}                                              		       
\renewcommand\vec[1]{\accentset{\,\rightarrow}{#1}}                              
\newcommand\spacevec[1]{\accentset{\,\rightarrow}{#1}}                        
\newcommand\contravec[1]{\tilde{ #1}}                                                     
\newcommand\contraspacevec[1]{\spacevec{\tilde{#1}}}                          
\newcommand\spacevecg[1]{\vv #1}                                                         
\newcommand\spacestatevec[1]{\vv{\spacevec{#1}}}                               
\newcommand\statevec[1]{\mathbf #1}                                                     
\newcommand\statevecGreek[1]{\boldsymbol #1}                                     
\newcommand\contrastatevec[1]{\tilde{\mathbf #1}}                                 
\newcommand\acclrvec[1]{\accentset{\,\leftrightarrow}{#1}}                      
\newcommand\bigstatevec[1]{\acclrvec{{\mathbf #1}}}                              
\newcommand\biggreekstatevec[1]{\acclrvec{\boldsymbol #1}}                
\newcommand\bigcontravec[1]{\acclrvec{\tilde{\mathbf{#1}}}}                   
\newcommand\biggreekcontravec[1]{\acclrvec{\tilde{\boldsymbol #1}}}    
\newcommand\vecNabla{\accentset{\,\rightarrow}\nabla}                         
\newcommand\vecNablaXi{\accentset{\,\rightarrow}\nabla_{\!\xi}}            
\newcommand\vecNablaX{\accentset{\,\rightarrow}\nabla_{\!x}}              
\newcommand\mmatrix[1]{\underbar{#1}}				
\newcommand\fiveMatrix[1]{\mathsf{ #1}}                          
\newcommand\fifteenMatrix[1]{\underline{\mathsf{ #1}}}    
\newcommand\matrixvec[1]{\mathcal #1}                           
\newcommand\bigmatrix[1]{\mathfrak #1}                          
\newcommand{\dmat}{\matrixvec{D}}     
\newcommand{\qmat}{\matrixvec{Q}}    
\newcommand{\mmat}{\matrixvec{M}}   
\newcommand{\bmat}{\matrixvec{B}}    
\newcommand\interiorfaces{\genfrac{}{}{0pt}{}{\mathrm{interior}}{\mathrm{faces}}}
\newcommand\boundaryfaces{\genfrac{}{}{0pt}{}{\mathrm{boundary}}{\mathrm{faces}}}
\newcommand\allfaces{\genfrac{}{}{0pt}{}{\mathrm{all}}{\mathrm{faces}}}
\newcommand{\testfuncOne}{\statevecGreek{\varphi}}  
\newcommand{\testfuncTwo}{\biggreekstatevec{\vartheta}}
\newcommand{\testfuncPhi}{\boldsymbol\phi}
%
\newcommand{\DD}{\spacevec{\mathbb{D}}\cdot}
\newcommand{\DDs}{\spacevec{\mathbb{D}}^{s}\cdot}
\newcommand{\IN}[1]{\mathbb I^{N}\!\!\left(#1\right)} 
\newcommand{\PN}[1]{\mathbb P^{#1}}
\newcommand{\LTwo}[1]{\mathbb L^{2}\!\left(#1\right)}
%
%
\newcommand\overRe{\frac{1}{{\operatorname{Re} }}}
\newcommand\twooverRe{\frac{2}{{\operatorname{Re} }}}
%
%
\newcommand\oneHalf{\frac{1}{2}}
\newcommand\oneFourth{\frac{1}{4}}

\begin{frontmatter}

\title{Global Bounds for the Error in Solutions of Linear Hyperbolic Systems due to Inaccurate Boundary Geometry}


\author[1,2]{David A. Kopriva}

\author[3]{Andrew R. Winters\corref{cor2}}
\ead{andrew.ross.winters@liu.se}

\author[3,4]{Jan Nordstr\"om}

\address[1]{Department of Mathematics, The Florida State University, Tallahassee, FL 32306, USA}
\address[2]{Computational Science Research Center, San Diego State University, San Diego, CA, USA}
\address[3]{Department of Mathematics, Applied Mathematics, Link\"{o}ping University, 581 83 Link\"{o}ping, Sweden}
\address[4]{Department of Mathematics and Applied Mathematics
 University of Johannesburg
 P.O. Box 524, Auckland Park 2006, South Africa.}

\cortext[cor2]{Corresponding author}




\begin{abstract}
\textcolor{black}{Meshes approximate the boundaries of a geometry when the boundaries are curved. The accuracy of the mesh then affects the error of computations of initial boundary value problems for partial differential equations, especially when using high order methods. Here, we} derive global estimates for the error in solutions of linear hyperbolic systems due to inaccurate boundary geometry. We show that the error is bounded by data and bounded in time when the solutions in the true and approximate domains are bounded. Just evaluating boundary data at the correct location has a secondary effect on the error, whereas the primary errors are from the Jacobian and metric terms. In two space dimensions, specifically, we show that to lowest order the errors are proportional to the errors in the boundary curve locations {\it and} their derivatives. \textcolor{black}{Therefore, high order accuracy computations cannot be obtained unless the mesh is also high order.} The results illustrate the importance of accurately approximating boundaries and should be helpful guides for high-order mesh generation for advection-dominated problems and the design of optimization algorithms for boundary approximations.
\end{abstract}

\end{frontmatter}

\section{Introduction}
The quest to compute more accurate solutions of advection-dominated problems in fluid mechanics, electrodynamics, etc. has led to the development of stable high-order methods, such as spectral element \cite{Canuto:2007fj} or high-order summation-by-parts finite difference methods, cf.~\cite{Canuto:2007fj,svard2014review}. Spectral element methods, in particular, are now mature, with a number of open source implementations available to solve complex physical problems in complex domains \cite{ampuero2012spectral,CANTWELL2015205,FERRER2023108700,FISCHER2022102982,KRAIS2021186,Martire:2021xt,KURZ2025109388,schlottkelakemper2021purely}.

The workflow to compute numerical approximations to the solution of initial boundary-value problems (IBVPs) for a partial differential equation with one of these methods starts with the definition of the geometry of an arbitrary domain, $\Omega$, which is then used by a mesh generator to generate a mesh. The mesh generator will approximate the geometry. For finite element or spectral element methods, boundary curves or surfaces in the mesh will usually be approximated by polynomials.  Spectral element methods will include large numbers of degrees of freedom corresponding to high (4-20+) polynomial orders. The mesh will then be passed to the solver, where each element will be mapped from a reference element, $\mathcal D$, using the approximate boundary information supplied by the mesher through a mapping $\spacevec X^e: \mathcal D \rightarrow \Omega_e$ that {\it approximates} the correct mapping, $\spacevec X: \mathcal D \rightarrow \Omega$. For finite difference schemes, the approximate mapping function would be defined to map a regular grid, $\mathcal D$ to the entire $\Omega_e$, or to a single block in a multiblock decomposition. We therefore view the domains $\Omega$ and $\Omega_e$ as representing an element, a block, or the whole domain.

The solver then computes its approximate solution on the reference domain with the approximate transformation, which we notate as $Q(\mathcal D;\spacevec X^e)$. $Q$ approximates the correct solution, $q(\mathcal D;\spacevec X)$,                           on the reference domain with the correct transformation.
The solver itself is only approximate and produces an error. We can formally write the total error of the approximation on an element or finite difference domain as
\begin{equation}
\begin{split}
    e &= Q\left(\mathcal D;\spacevec X^e\right) - q\left(\mathcal D;\spacevec X\right) 
    \\&= \left\{q\left(\mathcal D;\spacevec X^e\right) - q\left(\mathcal D;\spacevec X\right)\right\}
    +\left\{Q\left(\mathcal D;\spacevec X^e\right) - q\left(\mathcal D;\spacevec X^e\right)\right\}
    \\& = e_{mesh} + e_{solver}.
    \end{split}
    \label{eq:SplitError}
\end{equation}
The first term represents the {\it solution error} caused by the geometry error, called the ``mesh generation error". The second represents the error of the approximation and is due to the solver: finite difference, finite volume, finite element, spectral element, etc. The second term is what numerical convergence analysis studies, under the notion that $\spacevec X^e = \spacevec X$. Written as in \eqref{eq:SplitError}, we see that to get an accurate solution, one needs an accurate representation of the geometry through the mesh mapping {\it in addition to} an accurate solution on the mesh. 

We interpret \eqref{eq:SplitError} in two ways. 
\begin{enumerate}
    \item 
The total error is the sum of the mesh and the solver errors. 
Examples have shown that the mesh generation error, $e_{mesh}$, can be larger than $e_{solver}$ produced by the solvers themselves \cite{Hindenlang:2014gl,Minakowski:2020kq}, especially when using high order schemes. Ref. \cite{Minakowski:2020kq} in particular provides rationale for understanding and control of this error. For a given mesh, a converged approximate solution will still have an error, $e = e_{mesh}$.
\item It suggests that one should try to minimize the first term by minimizing some measure of error due to the geometry, especially \textcolor{black}{when it is to be used by} high order methods \cite{ISI:000306588600006}. 
\end{enumerate}

Numerous efforts have been made in the mesh generation community to optimize boundary approximations to minimize the mesh generation errors, effectively to minimize $\spacevec X^e - \spacevec X$ in some sense,  e.g. \cite{ISI:000306588600006,RUIZGIRONES2015122,ISI:000172856600011,TOULORGE2016361} to cite only four. Different approaches have been used to perform these optimizations.  In \cite{TOULORGE2016361}, three optimization procedures were presented using the Hausdorff distance, Fr\'{e}chet distance, and a Taylor series estimate of the geometric error in the boundary approximations, the latter of which minimizes a functional of the normal and tangent vectors.  The paper \cite{RUIZGIRONES2015122} minimized an area-based disparity measure. The authors of \cite{ISI:000172856600011} found optimal locations of interpolation nodes by minimizing a functional of the square of the distance between the true and approximate points, whereas the authors of \cite{ISI:000306588600006} showed that an optimization of the curve parameterization that ensured that the boundary normals were correct at the interpolation nodes gave orders of magnitude improvement in the solution error over other approaches. 

\textcolor{black}{None of those efforts address the central question:} How does the boundary approximation error, through $\spacevec X^e - \spacevec X$, created by the mesh, translate into the solution error, $q\left(\mathcal D;\spacevec X^e\right) - q\left(\mathcal D;\spacevec X\right)$, represented in the first term in \eqref{eq:SplitError}? Without such knowledge, it is not clear how to properly minimize $e_{mesh}$ by minimizing the geometry error.

Analysis of errors due to uncertainty or errors in the geometry has been studied for elliptic boundary-value problems for many years, e.g. \cite{Babuska:2002ix,Minakowski:2020kq,Xue:2005uo}. Such papers have also focused on low order elements which have low numbers of internal degrees of freedom.

For hyperbolic initial boundary-value problems, and for high order methods with large numbers of degrees of freedom with which to approximate boundaries, there have been virtually no analytic studies of how the grid generation error affects the solutions. Energy estimates for the influence of boundary errors for hyperbolic problems were derived in the short note \cite{NORDSTROM2016438}. The authors conclude that i) wavespeeds will be incorrect within the domain due to the fact that the eigenvalues of the system matrices will be different, ii) normals will be inaccurate, which can lead to incorrect specification of boundary conditions, iii) boundary conditions can be applied at the wrong points in space, and iv) boundary data may be applied at the wrong locations. Looking at the error in the divergence due to inaccurate boundary approximation, the authors of \cite{CHUN2022115261} concluded that errors in the boundary normals and the transformation Jacobian affect the divergence error.

In this paper, we go into detail and derive global estimates for the $e_{mesh}$ term in \eqref{eq:SplitError}  due to inaccurate boundaries for the solutions of linear hyperbolic systems. \textcolor{black}{The error, $e_{mesh}$, is associated with the continuous problem, and is {\it independent of the solver used}, see \eqref{eq:SplitError}}. We consider linear systems because only for those can an equation for the error be derived. \textcolor{black}{The results are directly applicable to problems in linear acoustics or electromagnetics in the solution of Maxwell's equations.} Alternatively, a common view is that it is sufficient to consider the properties of the constant coefficient problem as the frozen coefficient simplification of the linearization of a nonlinear problem \cite{GKOBook}.

We show that the error is bounded when the solutions in the correct and inaccurate domains are bounded. We also show that the primary solution errors are due to errors in the Jacobian and metric terms. Errors due to the evaluation of boundary functions at incorrect locations of the boundaries are secondary effects, being dependent on the product of small quantities. So it is not enough to simply evaluate boundary data at the correct points. In two space dimensions, specifically, we show that to lowest order the errors are proportional to the errors in the boundary curve locations {\it and} their derivatives. The results give guidance on how to optimize $\spacevec X^e - \spacevec X$ to minimize $e_{mesh}$ for the solutions of advection-dominated problems.

\section{Notation and Nomenclature}

\textcolor{black}{The derivations in this work are quite technical. For clarity and easy referral, we collect the notation that we use throughout the paper here. Quantities referring to the erroneous/approximated domain, $\Omega_e$, like Jacobian $J_e$, or mapping, $\spacevec X^e$, are denoted with a sub or superscript, $e$. Other quantities referring to the correct/true domain, $\Omega$, or the reference domain, $\mathcal D$, are denoted without sub or superscripts. Notation for vectors, matrices, and norms and other objects are also defined here. Finally, we list the various error quantities used.}

{\centering
\begin{longtable}[c]{ll}
$d$ & Number of space dimensions $\le 3$\\
$\Omega$ & Correct domain $\subset \mathcal R^d$\\
$\Omega_e$ & Erroneous domain $\subset \mathcal R^d$\\
$\mathcal D$ & Reference domain $= [0,1]^d$\\
$\spacevec x$ & Space vector, $\spacevec x \in \mathcal R^d$\\
$\contravec x$ & Contravariant space vector, $\contravec x \in \mathcal D$\\
$\statevec q$ & State vector\\
$\mmatrix A$ & Matrix \\
$\mmatrix A^\pm$ & Characteristic matrices, $\mmatrix A^\pm = \oneHalf\left(\mmatrix A \pm |\mmatrix A|\right)$ \\
$\mmatrix R$ & Reflection matrix\\
$\cdot_e$ & Quantity on the erroneous domain\\
$\spacevec \xi$& Reference space coordinate, e.g.  $\spacevec \xi = (\xi^1,\xi^2,\xi^3)=(\xi,\eta,\zeta)$\\
$\spacevec X\left(\spacevec \xi\right)$ & Correct mapping, $\spacevec X:\mathcal D\rightarrow \Omega$\\
$\spacevec X^e\left(\spacevec \xi\right)$ & Erroneous mapping, $\spacevec X^e:\mathcal D\rightarrow \Omega_e$\\
$\spacevec a_i$ & $i^{th}$ covariant basis vector, $\spacevec a_i = \spacevec X_{\xi^i}$\\
$J\spacevec a^i$ & $i^{th}$ volume weighted contravariant basis vector, $J\spacevec a^i = \spacevec a_j\times\spacevec a_k,\; i,j,k$ cyclic\\
$J$ & Transformation Jacobian $J=\spacevec a_i\cdot(\spacevec a_j\times\spacevec a_k),\; i,j,k$ cyclic\\
$\hat n$ & Unit normal to the boundary of the reference domain\\
$\iprod{\cdot,\cdot}$ & $\mathcal L^2$ inner product\\
$\inorm{\cdot}_J$ &Volume weighted $\mathcal L^2$ norm, $\iprod{J\cdot,\cdot}^\oneHalf$\\
$\statevec e$ & Solution error\\
$\rho$ & Ratio of correct and erroneous Jacobians, $J/J_e$\\
$\epsilon$ & Error in the Jacobian, $\epsilon = J_e - J$\\
$\varepsilon$ & Error in the Jacobian ratio, $\rho$. $\varepsilon = J/J_e - 1 = \rho - 1 = -\epsilon/(J + \epsilon)$\\
$N$ & Polynomial order\\
\end{longtable}
\setcounter{table}{0} 
}

\section{Energy Bounds on the Correct and Erroneous Domains}

We first introduce the notation we use here and the energy method to compute the energy bounds for a symmetric, constant-coefficient linear hyperbolic system
\begin{equation}
\statevec q_t + \nabla_x\cdot \left(\spacevec{\mmatrix A}\statevec q\right) = 0,
\end{equation}
where $\statevec q$ is the state vector of dimension $N_{eq}$, the number of equations, and
\begin{equation}
\spacevec{\mmatrix A} = \sum_{i=1}^d \mmatrix A_i\hat x_i = \mmatrix A_1 \hat x + \mmatrix A_2 \hat y + ...
\end{equation}
is a $d$-dimensional space vector of the coefficient matrices of size $N_{eq}\times N_{eq}$, with $\hat x_i$ being a unit vector in the $i^{th}$ coordinate direction. With this notation, product $\spacevec{\mmatrix A}\statevec q = \sum_{i=1}^d (\mmatrix A_i\statevec q)\hat x_i$ is a space vector of state vectors (a vector of fluxes, in this particular case). The dot product of $\spacevec{\mmatrix A}$ with a space vector, $\spacevec n= n_1\hat x + n_2\hat y+\ldots$ is a matrix, $\mmatrix A = \sum_{i=1}^d A_i n_i$, from which it follows that the divergence defined as 
\begin{equation}
    \nabla_x\cdot \left(\spacevec{\mmatrix A}\statevec q\right) = \sum_{i=1}^d \frac{\partial\left(\mmatrix A_i\statevec q\right)}{\partial x_i},
\end{equation}
is a state vector, see, e.g. \cite{BR12018}.

\begin{figure}[htbp] 
   \centering
   \includegraphics[width=0.6\textwidth]{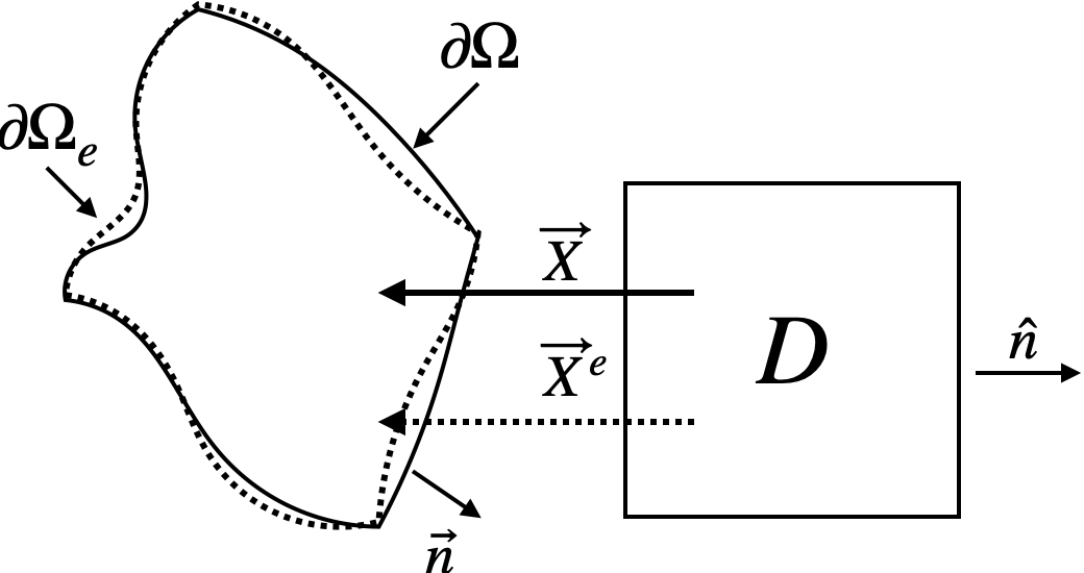} 
   \caption{Diagram of the correct, $\Omega$, erroneous, $\Omega_e$, and reference, $\mathcal D$, domains, \textcolor{black}{shown in two space dimensions}}
   \label{fig:CorrectAndErroneousMapping}
\end{figure}

The equations are posed on \textcolor{black}{an arbitrary domain $\Omega\subset \mathcal R^d$}, which we call the ``correct" domain, and so $\statevec q$ corresponds to the correct solution $q(\mathcal D;\spacevec X)$ described in the introduction. However the correct domain must almost always be approximated, especially when the boundaries are curved, which means that numerically one actually solves on the approximate, or ``erroneous" domain $\Omega_e$, see Fig. \ref{fig:CorrectAndErroneousMapping}. In the forthcoming analysis, we assume that both the correct and erroneous domains can be mapped from a reference domain $\spacevec \xi \in\mathcal D = [0,1]^d$ by a single {\it smooth} mapping $\spacevec x = \spacevec X\left(\spacevec \xi\right)$ for the correct domain and $\spacevec x = \spacevec X^e\left(\spacevec \xi\right)$ for the erroneous one. We require smoothness so that derivatives of the error quantities are continuous. 

\begin{remark}
As described in the introduction, the mapping can be considered to refer to a whole domain, as would be used in a mapped finite difference approximation, or just a single boundary element or boundary block in a mesh.\end{remark}

To complete the hyperbolic initial value problem (IBVP) we must introduce boundary and initial conditions. Here we consider two types of boundary conditions. The first sets external data, corresponding to inflow/outflow type conditions. The second sets pure reflection conditions, which correspond to wall-type boundaries. The initial conditions are specified on each domain by a function $\statevec q(\spacevec x,0) = \statevec q_0(\spacevec x)$.

\subsection{Energy Bounds for the Correct Domain}
From the transformation $\spacevec x = \spacevec X\left(\spacevec \xi\right)$, we can define the covariant basis and volume weighted contravariant basis vectors
\[
\spacevec a_i = \frac{\partial \spacevec X}{\partial \xi^i},\quad J\spacevec a^i = \spacevec a_j\times \spacevec a_k,\; i = 1,2,3,\; i,j,k \;cyclic,
\]
where, in three space dimensions, $\spacevec \xi = (\xi^1,\xi^2,\xi^3) \equiv (\xi,\eta,\zeta)$. The Jacobian of the transformation itself is $J = \spacevec a_i 
\cdot(\spacevec a_j \times \spacevec a_k)$, again with $i,j,k$ cyclic.

The equations for the correct domain transform onto the reference domain as
\begin{equation}
J \statevec q_t + \nabla_\xi\cdot \left(\contraspacevec{\mmatrix A}\statevec q\right) = 0,
\label{eq:CorrectStrong}
\end{equation}
where the matrix components of the space vector of the contravariant coefficient matrices are $\contraspacevec{\mmatrix A}^i = J\spacevec a^i\cdot\spacevec {\mmatrix A},\; i=1,2,3$, so that
\begin{equation}
    \contraspacevec{\mmatrix A} = \sum_{i=1}^d \contraspacevec{\mmatrix A}^i\hat \xi_i,
\end{equation}
where $\hat \xi_i$ is unit coordinate vector of the $i^{th}$ coordinate direction in reference space.
We assume that the correct domain is ``valid" in that $J_{max}\ge J\ge J_{min}>0$. Due to the metric identities, $\sum_{i=1}^d\frac{\partial(Ja^i_n)}{\partial\xi^i}=0,\;n=1,2,\ldots,d$, and the fact that $\spacevec{\mmatrix A}$ is constant, $\nabla_\xi\cdot \contraspacevec{\mmatrix A} = 0$, cf. \cite{kopriva2006}.

The solution energy is bounded by boundary and initial data. To show this, let $\statevecGreek\phi \in \mathcal L^2$ be a test function and let
\begin{equation}
\iprod{\statevec q,\statevecGreek\phi} = \int_\mathcal D \statevec q^T\statevecGreek\phi d\spacevec\xi
\end{equation}
be the $\mathcal L^2$ inner product.
Multiplying \eqref{eq:CorrectStrong} by $\statevecGreek\phi^T$ and integrating over the domain, it can be written in inner product notation as
\begin{equation}
\iprod{J\statevec q_t,\statevecGreek\phi} + \iprod{\nabla_\xi\cdot\left(\contraspacevec{\mmatrix A}\statevec q\right),\statevecGreek\phi} = 0.
\end{equation}
We then integrate the divergence term by parts to separate the boundary from the interior. Since the symmetric coefficient matrices are constant, and the divergence of the metric terms is zero,
\begin{equation}
\iprod{J\statevec q_t,\statevecGreek\phi} + \int_{\partial\mathcal D} \statevecGreek\phi^T \contraspacevec{\mmatrix A}\cdot\hat n\statevec q\dS 
- \iprod{\statevec q,\nabla_\xi\cdot\left(\contraspacevec{\mmatrix A}\statevecGreek\phi\right)} = 0.
\end{equation}

To obtain the energy, one replaces $\statevecGreek \phi\leftarrow \statevec q$ so that $\iprod{J\statevec q_t,\statevec q} = \oneHalf\ddt \inorm{\statevec q}^2_J$. Let the normal coefficient matrix be defined as $\mmatrix A \equiv \contraspacevec{\mmatrix A}\cdot\hat n$, where $\hat n$ is the outward normal to the reference domain. The volume term is converted to a surface term. With $\mmatrix A$ being constant and symmetric,
\begin{equation}
\iprod{\statevec q,\nabla_\xi\cdot\left(\contraspacevec{\mmatrix A}\statevec q\right)} = \int_{\partial\mathcal D} \statevec q^T \mmatrix A \statevec q\dS - \iprod{\nabla_\xi\cdot\left(\contraspacevec{\mmatrix A}\statevec q\right),\statevec q},
\end{equation}
becomes
\begin{equation}
\iprod{\statevec q,\nabla_\xi\cdot\left(\contraspacevec{\mmatrix A}\statevec q\right)} =\oneHalf \int_{\partial\mathcal D} \statevec q^T \mmatrix A \statevec q\dS.
\label{eq:VoltoSurfIntegral}
\end{equation}
Then the energy satisfies
\begin{equation}
 \frac{d}{dt}\inorm{\statevec q}^2_{J} +  \int_{\partial\mathcal D} \statevec q^T \mmatrix A\statevec q\dS  = 0.
\label{eq:EnergyBeforeBC}
\end{equation}

Since we incorporate the two types of boundary conditions, we group the boundary into two parts, $\partial D = \partial D_g \bigcup \partial D_r$, where external boundary data is applied along $\partial D_g$, and reflection conditions are applied along $\partial D_r$.

Along $\partial D_g$, we apply characteristic boundary conditions and split the normal coefficient matrix according to its positive and negative eigenvalues as $\mmatrix A = \mmatrix A^+ + \mmatrix A^-$ where $\mmatrix A^\pm = \oneHalf(\mmatrix A \pm |\mmatrix A|)$.  Since the coefficient matrices are symmetric, we can write
\begin{equation}
\mmatrix A = \mmatrix P \, \mmatrix{$\Lambda$} \, \mmatrix P^T = \mmatrix P \, \mmatrix{$\Lambda$}^+ \mmatrix P^T + \mmatrix P \, \mmatrix{$\Lambda$}^- \mmatrix P^T\equiv \mmatrix A^+ + \mmatrix A^- ,
\end{equation} 
where
\begin{equation}
\mmatrix{$\Lambda$}^+ = \left[\begin{array}{cc}\bar{\mmatrix{$\Lambda$}}^+ &\mmatrix 0 \\ \mmatrix0 &\mmatrix 0\end{array}\right], \quad \mmatrix{$\Lambda$}^- = \left[\begin{array}{cc}\mmatrix 0 &\mmatrix 0 \\ \mmatrix 0 &\bar{\mmatrix{$\Lambda$}}^-\end{array}\right].
\end{equation}
The matrices $\bar{\mmatrix{$\Lambda$}}^\pm$ contain the positive and negative eigenvalues of $\mmatrix A$. (We assume no zero eigenvalues, but if they are present they do not affect the results.)
Then
\begin{equation}
\int_{\partial\mathcal D_g} \statevec q^T \mmatrix A \statevec q\dS = \int_{\partial\mathcal D_g} \statevec q^T \mmatrix A^+\statevec q\dS - \int_{\partial\mathcal D_g} \statevec q^T |\mmatrix A^-| \statevec q\dS.
\label{eq:ExternalBSplit}
\end{equation}
Specifying external data is equivalent to replacing $\statevec q = \statevec g(\spacevec x,t)$ in the $|\mmatrix A^-|$ term in \eqref{eq:ExternalBSplit}, where $\statevec g$ is boundary data,
\begin{equation}
\int_{\partial\mathcal D_g} \statevec q^T \mmatrix A \statevec q\dS \leftarrow \int_{\partial\mathcal D_g} \statevec q^T \mmatrix A^+\statevec q\dS - \int_{\partial\mathcal D_g} \statevec g^T \left(\spacevec X\left(\spacevec \xi\right),t\right)|\mmatrix A^-| \statevec g\left(\spacevec X\left(\spacevec \xi\right),t\right)\dS,
\label{eq:BoundaryReplacement}
\end{equation}
see \cite{Nordstrom:2016jk}, \textcolor{black}{where the symbol ``$\leftarrow$" represents replacement. With the relation \eqref{eq:BoundaryReplacement}, the boundary terms are limited by external data.}

The second boundary condition along $\partial \mathcal{D}_r$ represents reflection, and includes wall-type conditions.
Here, characteristic boundary data is also used to determine the boundary terms, writing
\begin{equation}
\statevec q^T\mmatrix A\statevec q 
= \statevec q^T\mmatrix P \, \mmatrix{$\Lambda$} \left(\mmatrix P^T\statevec q\right) 
\equiv \statevec w^T \mmatrix{$\Lambda$}\statevec w 
=\left[\begin{array}{c} \statevec w^+ \\\statevec w^-\end{array}\right]^T  \left[\begin{array}{cc}\bar{\mmatrix{$\Lambda$}}^+ &\mmatrix 0 \\ \mmatrix 0 & \bar{\mmatrix{$\Lambda$}}^-\end{array}\right]\left[\begin{array}{c} \statevec w^+ \\\statevec w^-\end{array}\right],
\end{equation}
where $\statevec w= \mmatrix P^T\statevec q$ and $\statevec w^\pm$ are the characteristic variables grouped according to the positive and negative eigenvalues of $\mmatrix A$.

The general form of reflection boundary conditions without sources are of the form
\begin{equation}
\statevec w^- = \bar{\mmatrix R}\statevec w^+,
\label{eq:GenReflection}
\end{equation}
where $\bar{\mmatrix R}$ is the reflection matrix. Then, with the boundary condition, \eqref{eq:GenReflection}, the reflected state is
\begin{equation}
\statevec w_R = \left[\begin{array}{cc}\mmatrix I & \mmatrix 0 \\\bar{\mmatrix R} & \mmatrix 0\end{array}\right]\statevec w\equiv \mmatrix R\statevec w.
\label{eq:GenReflW}
\end{equation}
Inserting the boundary condition \eqref{eq:GenReflW} into the boundary integrand,
\begin{equation}
\statevec q^T \mmatrix A \statevec q \leftarrow \statevec w^T \mmatrix{$\Lambda$}^+\statevec w + 
\statevec w ^T \mmatrix R^T\mmatrix{$\Lambda$}^-\mmatrix R\statevec w
 = \statevec q^T\left[\mmatrix P\left(\mmatrix{$\Lambda$}^+ + \mmatrix R^T \mmatrix{$\Lambda$}^-\mmatrix R\right)\mmatrix P^T\right]\statevec q
\end{equation}
along reflection boundaries, so the boundary term is non-negative if the reflection matrix $\mmatrix R$ satisfies
\begin{equation}
\mmatrix{$\Lambda$}^+ + \mmatrix R^T \mmatrix{$\Lambda$}^-\mmatrix R\ge 0.
\label{eq:CorrectReflectionCondition}
\end{equation}
In this paper we will assume that the reflection matrix is specified so that \eqref{eq:CorrectReflectionCondition} is true.
\begin{remark}
    Setting $\mmatrix R=0$ is equivalent to setting trivial boundary data, $\statevec g = 0$.
\end{remark}

Gathering the results,
\begin{equation}
\begin{split}
 \frac{d}{dt}\inorm{\statevec q}^2_{J} &+  \int_{\partial\mathcal D_r}  \statevec q^T\left[\mmatrix P\left(\mmatrix{$\Lambda$}^+ + \mmatrix R^T \mmatrix{$\Lambda$}^-\mmatrix R\right)\mmatrix P^T\right]\statevec q\dS +\int_{\partial\mathcal D_g} \statevec q^T \mmatrix A^+\statevec q\dS 
 \\& =\int_{\partial\mathcal D_g} \statevec g^T \left(\spacevec X\left(\spacevec \xi\right),t\right)|\mmatrix A^-| \statevec g\left(\spacevec X\left(\spacevec \xi\right),t\right)\dS .
 \end{split}
\label{eq:EnergyAfterReflectionBC}
\end{equation}
The boundary terms on the left of \eqref{eq:EnergyAfterReflectionBC} are dissipative, so, after we integrate in time and apply the initial condition $\statevec q(\spacevec x,0) = \statevec q_0(\spacevec x)$,
\begin{equation}
\inorm{\statevec q(t)}^2_{J}\le \inorm{\statevec q_0\left(\spacevec X\left(\spacevec \xi\right)\right)}^2_{J} + \int_0^t\int_{\partial\mathcal D_g} \statevec g^T \left(\spacevec X\left(\spacevec \xi\right),t\right)|\mmatrix A^-| \statevec g\left(\spacevec X\left(\spacevec \xi\right),t\right)\dS dt.
\label{eq:CorrectSolutionBound}
\end{equation}
Thus, the energy is bounded by initial and boundary data.
\begin{remark}
The case 
\begin{equation}
\mmatrix{$\Lambda$}^+ + \mmatrix R^T \mmatrix{$\Lambda$}^-\mmatrix R =0
\label{eq:PerfectReflectionCorrect}
\end{equation}
produces no loss or growth of energy. For instance if $\partial D_g = \emptyset$, then  
\begin{equation}
\inorm{\statevec q(t)}_{J}= \inorm{\statevec q_0\left(\spacevec X\left(\spacevec \xi\right)\right)}_{J},
\end{equation}
so \eqref{eq:PerfectReflectionCorrect} corresponds to a perfectly reflecting boundary condition.
\end{remark}

\subsection{Energy Bounds for the Erroneous Domain}

Now we assume that errors are introduced along the boundary of the correct domain, e.g. due to approximating boundary curves and surfaces, creating the ``erroneous" domain, $\Omega_e$. In that case, there is a different mapping, $\spacevec x = \spacevec X^e\left(\spacevec \xi\right)$ from the reference domain to the erroneous domain.  We assume that the erroneous domain is still valid in that $J_{e,max}\ge J_e \ge J_{e,min} >0$. For the erroneous problem, let $\contraspacevec{\mmatrix A}_e^{i} = J\spacevec a_e^i\cdot\spacevec {\mmatrix A}$ so that the erroneous solution, $\statevec v$, which we use to represent the generic $q(\mathcal D;\spacevec X^e)$, satisfies
\begin{equation}
J_e \statevec v_t + \nabla_\xi\cdot \left(\contraspacevec{\mmatrix A}_e\statevec v\right) = 0.
\label{eq:ErroneousStrong1}
\end{equation}
It follows also that  $\nabla_\xi\cdot \contraspacevec{\mmatrix A}_e = 0$. 

Following the same steps as for the solution on the correct domain, the solution on the erroneous domain is also bounded by data if boundary data is specified and exists,
\begin{equation}
\inorm{\statevec v(t)}^2_{J_e}\le \inorm{\statevec q_0\left(\spacevec X^e\left(\spacevec \xi\right)\right)}^2_{J_e} + \int_0^t  \int_{\partial\mathcal D_g} \statevec g^T\left(\spacevec X^e\left(\spacevec \xi\right),t\right) |\mmatrix A_e^-| \statevec g\left(\spacevec X^e\left(\spacevec \xi\right),t\right)\dS,
\label{eq:ErroneousEnergyBound}
\end{equation}
provided that any reflection matrix $\mmatrix R_e$ satisfies
\begin{equation}
\mmatrix{$\Lambda$}_e^+ + \mmatrix R_e^T \mmatrix{$\Lambda$}_e^-\mmatrix R_e\ge 0.
\label{eq:ErroneousREflectionCondition}
\end{equation}

Several technical considerations should be noted for \eqref{eq:ErroneousEnergyBound} to be valid. First, the boundary values and initial conditions must exist on both domains, that is,
the bounds require that the boundary conditions are defined on both $\partial \Omega$ and $\partial\Omega_e$, and initial conditions are defined on both $\Omega$ and $\Omega_e$. This can be either trivial, e.g. constant inflow conditions, or not, and might require some sort of extension of the conditions. The extension itself may be completely erroneous.

\begin{example}
\textcolor{black}{The boundary data $\statevec g$ may not even be defined for some values of $\spacevec X^e$, leading to no bound in \eqref{eq:ErroneousEnergyBound}. Suppose $\min_x \Omega = 0$ and the scalar boundary data is $\statevec g(x) = \sqrt{x}$. Then for the erroneous domain, $\statevec g(\spacevec{X}^e) = \sqrt{{X}^e}$. If $\partial \Omega_e$ is approximated as a high order polynomial interpolation of $\partial\Omega$, oscillations in the polynomial will produce sections where $X^e<0$, and $\statevec g$ will not exist at those points, cf. Fig. \ref{fig:CorrectAndErroneousMapping}. The same holds true with regards to the initial condition, $\statevec q_0$.}
\end{example}

Next, the boundary segmentations, $\partial D_g$ and $\partial D_r$ do not have to be the same, since the parametrizations can differ so that at a particular value of $\spacevec \xi$ in one domain may have specified conditions while the other has reflections. To simplify the analysis, we require that the boundary segmentations be the same for both domains. In practice, this is not a concern since domains (and elements) are usually defined so that only one boundary condition is applied to any parametrized segment of the boundary, e.g. one side of an element. 

Finally, the energy may also not be bounded if the exact reflection conditions $\mmatrix R$ are applied on the erroneous boundary. For $\partial D_g = \emptyset$ with $\mmatrix R_e = \mmatrix R$, the energy method produces
\begin{equation}
 \frac{d}{dt}\inorm{\statevec v}^2_{J} +  \int_{\partial\mathcal D}  \statevec v^T\left[\mmatrix P_e\left(\mmatrix{$\Lambda$}_e^+ + \mmatrix R^T \mmatrix{$\Lambda$}_e^-\mmatrix R\right)\mmatrix P_e^T\right]\statevec v\dS  = 0,
\label{eq:ErroneousEnergyAfterReflectionBC}
\end{equation}
{\it but only if the number of positive and negative eigenvalues do not differ between the correct and erroneous domains}. (These can differ depending on the differences in the normals.) If they do, then the reflection matrix doesn't have the correct dimensions and the problem is ill-posed. Even so,
if we write 
$
\mmatrix{$\Lambda$}^+_e = \mmatrix{$\Lambda$}^+ + \Delta \mmatrix{$\Lambda$}^+
$
, etc., then
\begin{equation}
\mmatrix{$\Lambda$}_e^+ + \mmatrix R^T \mmatrix{$\Lambda$}_e^-\mmatrix R = \left(\mmatrix{$\Lambda$}^+ + \mmatrix R^T \mmatrix{$\Lambda$}^-\mmatrix R\right) + \left(\Delta\mmatrix{$\Lambda$}^+ + \mmatrix R^T \Delta\mmatrix{$\Lambda$}^-\mmatrix R \right).
\end{equation}
Therefore the energy in \eqref{eq:ErroneousEnergyAfterReflectionBC} is only bounded if
\begin{equation}
\left(\Delta\mmatrix{$\Lambda$}^+ + \mmatrix R^T \Delta\mmatrix{$\Lambda$}^-\mmatrix R \right) +\left(\mmatrix{$\Lambda$}^+ + \mmatrix R^T \mmatrix{$\Lambda$}^-\mmatrix R\right)\ge 0.
\end{equation}
If, for instance, perfect reflection conditions are specified for the correct domain, then boundedness requires that
\begin{equation}
\left(\Delta\mmatrix{$\Lambda$}^+ + \mmatrix R^T \Delta\mmatrix{$\Lambda$}^-\mmatrix R \right) \ge 0,
\end{equation}
which, depending on the boundary errors, does not have to hold. To ensure that the erroneous problem is energy bounded, we require that the reflection matrix, $R_e$, satisfies \eqref{eq:ErroneousREflectionCondition}. In practice, i.e. in numerical simulations, a reflection condition is usually implemented using the available (erroneous) geometry, thus ensuring \eqref{eq:ErroneousREflectionCondition}.

\begin{remark}
    We see, then, that a well-posed problem on the correct domain does not {\rm necessarily} imply well-posedness on the erroneous domain. But we will not consider these ``worst case" scenarios here. Instead, we require the more likely situations where conditions on the initial and boundary values are defined so that the problems on both domains are well-posed.
\end{remark}

To ensure that the error is a continuous function of the mapping,  we will also require that the boundary and initial data are smooth. For example, the initial condition $\statevec q_0(\spacevec x)$ is defined on both the correct and erroneous domains. Then by the mean value theorem,
\begin{equation}
\statevec e(\spacevec \xi,0) \equiv \statevec q_0\left(\spacevec X^e\left(\spacevec\xi\right)\right) - \statevec q_0\left(\spacevec X\left(\spacevec\xi\right)\right) = \nabla_x \statevec q_0(\spacevec\nu)\cdot\left(\spacevec X^e\left(\spacevec\xi\right) - \spacevec X\left(\spacevec\xi\right)\right) 
\end{equation}
for some $\spacevec \nu$. For the initial error to depend continuously on the error in the mappings, the initial condition must be Lipschitz continuous. 
So, in addition to existence of boundary and initial data on both domains, we will require in the following that the error in the data is continuous with respect to the domain error and that all data have continuous and bounded first derivatives.

\section{The Error Equation}

The goal is to find how the erroneous solution, $\statevec v$, and the correct solution, $\statevec q$, differ. It is not easy to decide how to measure that error. (C.f the discussion in \cite{Minakowski:2020kq}.) If one attempts to measure it in physical space, then it is {\it undefined} in the regions where the two domains do not overlap. (See Fig. \ref{fig:CorrectAndErroneousMapping}.) On the other hand, the error is defined at all points in reference space. Although defining the error in reference space has the disadvantage that the difference between two values in reference space corresponds to the error at different physical space locations except in the limit as $\Omega_e\rightarrow\Omega$, it continuously converges to zero in that limit. For these reasons, we measure the error in reference space.

To that end, we multiply \eqref{eq:ErroneousStrong1} by the ratio of the Jacobians, $\rho \equiv J/J_e >0$, to get
\begin{equation}
J \statevec v_t + \rho\nabla_\xi\cdot \left(\contraspacevec{\mmatrix A}_e\statevec v\right) = 0
\label{eq:ErroneousStrongWeighted}
\end{equation}
for the erroneous solution.
We define the error equation by the difference between \eqref{eq:ErroneousStrongWeighted} and \eqref{eq:CorrectStrong},
\begin{equation}
J( \statevec v -\statevec q)_t  + \rho\nabla_\xi\cdot\left(\contraspacevec{\mmatrix A}_e\statevec v\right) - \nabla_\xi\cdot\left(\contraspacevec{\mmatrix A}\statevec q\right)= 0.
\end{equation}
For any $\mathcal L^2\cap C^0$ test function, $\statevecGreek\phi$, remembering that we are requiring extra smoothness in the data,
\begin{equation}
\iprod{J( \statevec v -\statevec q)_t,\statevecGreek\phi}  + \iprod{\rho\nabla_\xi\cdot\left(\contraspacevec{\mmatrix A}_e\statevec v\right),\statevecGreek\phi} - \iprod{\nabla_\xi\cdot\left(\contraspacevec{\mmatrix A}\statevec q\right),\statevecGreek\phi}= 0.
\label{eq:ErrorEquationOriginal}
\end{equation}

To separate erroneous from correct quantities, we define the error matrix $ \contraspacevec{\mmatrix E}$ through
\begin{equation}
\contraspacevec{\mmatrix A}_e = \contraspacevec{\mmatrix A} + \contraspacevec{\mmatrix E}.
\end{equation}
Since the coefficient matrix of the original problem is constant, the error matrix
\begin{equation}
\contraspacevec{\mmatrix E} = \sum_{i=1}^d \left\{J\spacevec a^i_e - J\spacevec a^i \right\}\cdot\spacevec{ \mmatrix A}\hat \xi^i
\label{eq:ErrorMatrices}
\end{equation}
depends only on the error in the metric terms. Note in particular that $\nabla_\xi\cdot \contraspacevec{\mmatrix E} = 0$ since in both cases the volume weighted contravariant basis vectors satisfy the metric identities.

 Let us also define the Jacobian error, $\varepsilon$, through $\rho = 1 + \varepsilon = 1 + \left(\frac{J}{J_e} - 1\right)$. Note that for valid domains, $\rho > 0$, so $\varepsilon > -1$.
In terms of the error matrix and Jacobian errors,
\begin{equation}
\begin{split}
 \iprod{\rho\nabla_\xi\cdot\left(\contraspacevec{\mmatrix A}_e\statevec v\right),\statevecGreek\phi} &= \iprod{ (1+\varepsilon)\nabla_\xi\cdot\left( \contraspacevec{\mmatrix A} + \contraspacevec{\mmatrix E}\right)\statevec v,\statevecGreek\phi}
 \\&=  \iprod{\nabla_\xi\cdot\left(\contraspacevec{\mmatrix A}\statevec v\right),\statevecGreek\phi} + \iprod{\varepsilon\nabla_\xi\cdot\left(\contraspacevec{\mmatrix A}\statevec v\right),\statevecGreek\phi} +  \iprod{ (1+\varepsilon)\nabla_\xi\cdot\left(\contraspacevec{\mmatrix E}\statevec v\right),\statevecGreek\phi}.
 \end{split}
 \label{eq:RewriteBMatrix}
\end{equation}
\begin{remark}Note that the second and third terms on the right of \eqref{eq:RewriteBMatrix} are expected to be small for small deformations of the domain, and the product of $\varepsilon$ and $\contraspacevec{\mmatrix E}$ will be secondary since $\varepsilon\nabla_\xi \cdot\left(\contraspacevec{\mmatrix E}\statevec v\right) = \left(\varepsilon \contraspacevec{\mmatrix E}\right)\cdot\nabla_\xi\statevec v$, since it includes the product of two small quantities.
\end{remark}

Substituting the formulation \eqref{eq:RewriteBMatrix}, the error equation \eqref{eq:ErrorEquationOriginal} becomes
\begin{equation}
\iprod{J( \statevec v -\statevec q)_t,\statevecGreek\phi} 
+ \iprod{\nabla_\xi\cdot\left(\contraspacevec{\mmatrix A}\left(\statevec v - \statevec q\right)\right),\statevecGreek\phi} 
+ \iprod{(1+\varepsilon)\nabla_\xi\cdot\left(\contraspacevec{\mmatrix E}\statevec v\right),\statevecGreek\phi} + \iprod{\varepsilon\nabla_\xi\cdot\left(\contraspacevec{\mmatrix A}\statevec v\right),\statevecGreek\phi}= 0.
\end{equation}
Let the error be $\statevec e \equiv \statevec v - \statevec q$, {\it which, again, is defined on the reference domain}. Then $\statevec v = \statevec q + \statevec e$ and
\begin{equation}
\begin{split}
\iprod{J \statevec e_t,\statevecGreek\phi} 
&+ \iprod{\nabla_\xi\cdot\left(\contraspacevec{\mmatrix A}\statevec e\right),\statevecGreek\phi} 
+ \iprod{(1+\varepsilon)\nabla_\xi\cdot\left(\contraspacevec{\mmatrix E}\statevec e\right),\statevecGreek\phi} + \iprod{\varepsilon\nabla_\xi\cdot\left(\contraspacevec{\mmatrix A}\statevec e\right),\statevecGreek\phi}
\\&+ \iprod{(1+\varepsilon)\nabla_\xi\cdot\left(\contraspacevec{\mmatrix E}\statevec q\right),\statevecGreek\phi} + \iprod{\varepsilon\nabla_\xi\cdot\left(\contraspacevec{\mmatrix A}\statevec q\right),\statevecGreek\phi}
= 0
\end{split}
\end{equation}
is the error equation.
When we set $\statevecGreek\phi = \statevec e$, we get an equation for the energy of the error,
\begin{equation}
\begin{split}
\oneHalf\frac{d}{dt}\inorm{\statevec e}^2_{J}
&+ \iprod{\nabla_\xi\cdot\left(\contraspacevec{\mmatrix A}\statevec e\right),\statevec e} 
+ \iprod{(1+\varepsilon)\nabla_\xi\cdot\left(\contraspacevec{\mmatrix E}\statevec e\right),\statevec e} + \iprod{\varepsilon\nabla_\xi\cdot\left(\contraspacevec{\mmatrix A}\statevec e\right),\statevec e}
\\&+ \iprod{(1+\varepsilon)\nabla_\xi\cdot\left(\contraspacevec{\mmatrix E}\statevec q\right),\statevec e} + \iprod{\varepsilon\nabla_\xi\cdot\left(\contraspacevec{\mmatrix A}\statevec q\right),\statevec e}
= 0.
\end{split}
\end{equation}

Using integration by parts, c.f. \eqref{eq:VoltoSurfIntegral},
\begin{equation}
\iprod{\nabla_\xi\cdot\left(\contraspacevec{\mmatrix A}\statevec e\right),\statevec e} = \oneHalf\int_{\partial\mathcal D} \statevec e^T \left(\contraspacevec{\mmatrix A}\cdot\hat n\right)\statevec e\dS,
\end{equation}
where, again, $\hat n$ is the unit outward normal on the reference domain.
Next,
\begin{equation}
\iprod{\varepsilon\nabla_\xi\cdot\left(\contraspacevec{\mmatrix A}\statevec e\right),\statevec e} = \oneHalf \int_{\partial\mathcal D} \varepsilon\statevec e^T\left(\contraspacevec{\mmatrix A}\cdot\hat n\right) \statevec e \dS -\oneHalf \iprod{\statevec e,\left(\contraspacevec{\mmatrix A}\cdot\nabla_\xi\varepsilon\right)\statevec e},
\label{eq:WeightedErrorInnerProduct}
\end{equation}
so this is like the previous term but with an extra volume term due to the gradient of $\varepsilon$.
Similarly,
\begin{equation}
 \iprod{\rho\nabla_\xi\cdot\left(\contraspacevec{\mmatrix E}\statevec e\right),\statevec e} = \oneHalf \int_{\partial\mathcal D} \rho\statevec e^T\left(\contraspacevec{\mmatrix E}\cdot\hat n\right) \statevec e \dS -\oneHalf \iprod{\statevec e,\left(\contraspacevec{\mmatrix E}\cdot\nabla_\xi\rho\right)\statevec e}.
\end{equation}
We put everything so far together, use the fact that $\rho = 1 + \varepsilon$, $\nabla_\xi\rho = \nabla_\xi\varepsilon$, and gather the boundary terms, to get
\begin{equation}
\begin{split}
\oneHalf\frac{d}{dt}\inorm{\statevec e}^2_{J}
&+ \oneHalf\int_{\partial\mathcal D} \statevec (1+\varepsilon)\statevec e^T
 \left\{\left(\contraspacevec{\mmatrix A}\cdot\hat n\right)  
+    \left(\contraspacevec{\mmatrix E}\cdot\hat n\right)  
 \right\} \statevec e \dS =
 \\&\oneHalf \iprod{\statevec e,\left(\contraspacevec{\mmatrix E}\cdot\nabla_\xi\varepsilon\right)\statevec e} +\oneHalf \iprod{\statevec e,\left(\contraspacevec{\mmatrix A}\cdot\nabla_\xi\varepsilon\right)\statevec e}
- \iprod{(1+\varepsilon)\nabla_\xi\cdot\left(\contraspacevec{\mmatrix E}\statevec q\right),\statevec e} - \iprod{\varepsilon\nabla_\xi\cdot\left(\contraspacevec{\mmatrix A}\statevec q\right),\statevec e}.
\end{split}
\label{eq:ErrorEnergy1}
\end{equation}

To get a bound on the error, boundary conditions must be applied to \eqref{eq:ErrorEnergy1}. Again, we will consider both externally specified data and reflection boundary conditions.

For externally specified data, the error on $\partial D_g$,
\[
\int_{\partial\mathcal D} \statevec (1+\varepsilon)\statevec e^T
 \left\{\left(\contraspacevec{\mmatrix A}\cdot\hat n\right)  
+    \left(\contraspacevec{\mmatrix E}\cdot\hat n\right)  
 \right\} \statevec e \dS = \int_{\partial\mathcal D} \statevec (1+\varepsilon)\statevec e^T
 \mmatrix A_e \statevec e \dS,
\]
where $\mmatrix A_e = \contraspacevec{\mmatrix A}_e\cdot\hat n = \contraspacevec{\mmatrix A}\cdot\hat n 
+   \contraspacevec{\mmatrix E}\cdot\hat n \equiv \mmatrix A + \mmatrix E$.
To specify characteristic boundary conditions with external data, one splits the matrix $\mmatrix A_e$ according to its eigenvalues to represent outgoing and incoming waves, for which the latter is specified by data \cite{Nordstrom:2016jk}. 
We will use the splitting $\statevec e^T \left(\mmatrix A^+ + \mmatrix E^+\right) \statevec e + \statevec e_g^T\left(\mmatrix A^-  + \mmatrix E^-\right)\statevec e_g$, since we can use it to relate $\mmatrix E^\pm$ to the boundary errors.

\textcolor{black}{To apply the boundary conditions, we replace the boundary integral by one split into incoming and outgoing waves, where incoming waves are specified by data,}
\begin{equation}
\begin{split}
\int_{\partial\mathcal D_g} (1 + \varepsilon)\statevec e^T
\mmatrix A_e   
 \statevec e \dS 
 \leftarrow \int_{\partial\mathcal D_g} (1 + \varepsilon)\statevec e^T
\left(\mmatrix A^+ + \mmatrix E^+\right)    
 \statevec e \dS - \int_{\partial\mathcal D_g} (1 + \varepsilon)\statevec e_g^T
(|\mmatrix A^-| + |\mmatrix E^- |)
 \statevec e_g \dS.
 \end{split}
 \label{eq:Boundary error splitting}
\end{equation}

Specified boundary data is given by external states $\statevec g =\statevec g\left(\spacevec x,t\right)$ applied to the incoming waves. We have two ways to specify the boundary data on the erroneous domain leading to the boundary error $\statevec e_g$:
\begin{itemize}
\item Evaluate the boundary data at the erroneous location, $\statevec g\left(\spacevec X^e\left(\spacevec \xi\right),t\right)$.

This is what one must do when the mesh generator, e.g. HOHQMesh \cite{Kopriva2024}, HOPR \cite{Hindenlang2015},  provides only an approximation to the boundary. This situation is where a boundary is defined by an interpolant only in terms of a set of points, and the exact ``truth'' boundary is never actually known by a solver.

\item Evaluate the data at the correct location, $\statevec g\left(\spacevec X\left(\spacevec \xi\right),t\right)$.

To evaluate at the correct location, the exact boundary must be supplied along with the mesh. 
In this case the boundary error, $\statevec e_g$, is zero and the last term in 
\eqref{eq:Boundary error splitting} vanishes.
\end{itemize}


In the first case, error is injected at the boundary of the reference domain
\begin{equation}
\statevec e_g = \statevec g\left(\spacevec X^e\left(\spacevec \xi\right),t\right) - \statevec g\left(\spacevec X\left(\spacevec \xi\right),t\right)\quad \spacevec \xi \in \partial D_g.
\end{equation}
For smooth boundary data, the mean value theorem allows the error to be written in terms of the boundary location error $\Delta \spacevec X$ as
\begin{equation}
\statevec e_g = \nabla_x\statevec g\left(\spacevec \nu,t\right)\cdot \left(\spacevec X^e\left(\spacevec \xi\right)
 - \spacevec X\left(\spacevec \xi\right)\right) =  \nabla_x\statevec g\cdot\Delta \spacevec X\left(\spacevec \xi\right),
 \label{eq:BoundaryErrorTerm}
\end{equation}
for some $\spacevec \nu$. 
Therefore, characteristic boundary terms in \eqref{eq:Boundary error splitting} become
\begin{equation}
\begin{split}
\int_{\partial\mathcal D_g}(1+\varepsilon) \statevec e^T
\mmatrix A_e   
 \statevec e \dS 
 \leftarrow &\int_{\partial\mathcal D_g} (1+\varepsilon)\statevec e^T
\left(\mmatrix A^+ + \mmatrix E^+\right)    
 \statevec e \dS 
 \\&- \int_{\partial\mathcal D_g} (1+\varepsilon)\Delta \spacevec X^T\cdot  \nabla_x\statevec g^T(|\mmatrix A^-| + |\mmatrix E^- |)  \nabla_x\statevec g\cdot\Delta \spacevec X \dS.
 \end{split}
 \label{eq:ErrorInjection}
\end{equation}
%

For reflection boundaries, we write the reflected error in the boundary term as
\begin{equation}
\int_{\partial\mathcal D_r} (1+\varepsilon)\statevec e^T 
\mmatrix A_e  \statevec e \dS
=
\int_{\partial\mathcal D_r}(1+\varepsilon) \statevec e^T
\left(\mmatrix A_e^+ +
\mmatrix A_e^- \right) \statevec e_r \dS.
\end{equation}
Let us define $\statevec W = \mmatrix P^T_e\statevec e$ to be the characteristic decomposition of the error with respect to the erroneous domain. Then
\begin{equation}
\int_{\partial\mathcal D_r}(1+\varepsilon) \left[\statevec e^T
\mmatrix A_e^+  \statevec e +  \statevec e_r^T\mmatrix A_e^-  \statevec e_r \right]\dS 
=\int_{\partial\mathcal D_r}(1+\varepsilon) \left[\statevec W^T\mmatrix {$\Lambda$}_e^+  \statevec W +  \statevec W^T\mmatrix {$\Lambda$}_e^-  \statevec W \right]\dS.
\end{equation}
To get the reflected error, we use the reflection matrix with respect to the erroneous domain and replace the characteristic variables in the $\mmatrix{$\Lambda$}^-_e$ term with
\begin{equation}
\statevec W_r = \mmatrix R_e\statevec W = \mmatrix R_e\mmatrix P_e^T\statevec e.
\end{equation}
So
\begin{equation}
\begin{split}
\int_{\partial\mathcal D_r}(1+\varepsilon) \left[\statevec W^T\mmatrix {$\Lambda$}_e^+  \statevec W +  \statevec W_r^T\mmatrix {$\Lambda$}_e^-  \statevec W_r \right]\dS 
&= \int_{\partial\mathcal D_r}(1+\varepsilon) \left[\statevec W^T\mmatrix {$\Lambda$}_e^+  \statevec W +  \statevec W^T\mmatrix R^T_e\mmatrix {$\Lambda$}_e^-  \mmatrix R_e\statevec W \right]\dS
\\&=\int_{\partial\mathcal D_r}(1+\varepsilon) \statevec W^T\left[\mmatrix {$\Lambda$}_e^+  + \mmatrix R^T_e\mmatrix {$\Lambda$}_e^-  \mmatrix R_e \right]\statevec W\dS\ge 0,
\end{split}
\end{equation}
when the reflection matrix on the erroneous domain is chosen properly as in \eqref{eq:ErroneousREflectionCondition}.

When the boundary data is evaluated at the correct point, the error injected along the boundary is zero, so the last term in \eqref{eq:ErrorInjection} vanishes. To distinguish between the two cases, we introduce a parameter $\gamma$ with
\begin{equation}
\gamma = \left\{  
\begin{gathered}
1\quad \mathrm{ evaluate \;along \;erroneous \;boundary} \hfill\\
0\quad \mathrm{ evaluate \;along \;correct \;boundary}\hfill
\end{gathered}\right. .
\label{eq:GammaDef}
\end{equation}
Then accounting for the boundary errors, we have the final form of the error energy equation,
\begin{equation}
\begin{split}
\oneHalf\frac{d}{dt}\inorm{\statevec e}^2_{J}
&+ \oneHalf\int_{\partial\mathcal D_g} (1+\varepsilon) \statevec e^T \left(\mmatrix A^+ + \mmatrix E^+\right)\statevec e \dS  +\oneHalf\int_{\partial\mathcal D_r}(1+\varepsilon) \statevec e^T\mmatrix P_e^T\left[\mmatrix {$\Lambda$}_e^+  + \mmatrix R^T_e\mmatrix {$\Lambda$}_e^-  \mmatrix R_e \right]\mmatrix P_e\statevec e\dS
\\&
=\frac{\gamma}{2}\int_{\partial\mathcal D_g} (1+\varepsilon)\Delta \spacevec X^T\cdot\nabla_x\statevec g^T  (|\mmatrix A^-| + |\mmatrix E^- |)  \nabla_x\statevec g\cdot\Delta \spacevec X \dS
 \\&
 +\oneHalf \iprod{\statevec e,\left(\contraspacevec{\mmatrix E}\cdot\nabla_\xi\varepsilon\right)\statevec e} 
 +\oneHalf \iprod{\statevec e,\left(\contraspacevec{\mmatrix A}\cdot\nabla_\xi\varepsilon\right)\statevec e}
\\&
- \iprod{(1+\varepsilon)\nabla_\xi\cdot\left(\contraspacevec{\mmatrix E}\statevec q\right),\statevec e} - \iprod{\varepsilon\nabla_\xi\cdot\left(\contraspacevec{\mmatrix A}\statevec q\right),\statevec e}.
\end{split}
\label{eq:CharSplitEnergyWGamma}
\end{equation}

To keep track of the numerous terms, we write \eqref{eq:CharSplitEnergyWGamma} in terms of the dissipation, $D$, boundary errors, $BG$, and volume errors, $V_1,\ldots,V_4$, as
\begin{equation}
    \inorm{\statevec e}_{J}\frac{d}{dt}\inorm{\statevec e}_{J} + D = \gamma BG + V1 + V2 + V3 + V4,
\label{eq:SimplifiedForm}
\end{equation}
where
\begin{equation}
\begin{gathered}
D=\oneHalf\int_{\partial\mathcal D_g}  (1+\varepsilon)\statevec e^T \mmatrix A_e^+ \statevec e \dS 
+ \oneHalf\int_{\partial\mathcal D_r}(1+\varepsilon) \statevec e^T\mmatrix P_e^T\left[\mmatrix {$\Lambda$}_e^+  + \mmatrix R^T_e\mmatrix {$\Lambda$}_e^-  \mmatrix R_e \right]\mmatrix P_e\statevec e\dS\ge 0,\hfill\\
BG = \frac{1}{2}\int_{\partial\mathcal D_g} (1+\varepsilon)\Delta \spacevec X^T\cdot\nabla_x\statevec g^T  (|\mmatrix A^-| + |\mmatrix E^- |)  \nabla_x\statevec g\cdot\Delta \spacevec X \dS,\hfill\\
V1 = \oneHalf \iprod{\statevec e,\left(\contraspacevec{\mmatrix E}\cdot\nabla_\xi\varepsilon\right)\statevec e},\quad
V2 = \oneHalf \iprod{\statevec e,\left(\contraspacevec{\mmatrix A}\cdot\nabla_\xi\varepsilon\right)\statevec e},\hfill\\
V3 = - \iprod{(1+\varepsilon)\nabla_\xi\cdot\left(\contraspacevec{\mmatrix E}\statevec q\right),\statevec e},\quad
V4 = - \iprod{\varepsilon\nabla_\xi\cdot\left(\contraspacevec{\mmatrix A}\statevec q\right),\statevec e}.\hfill
\end{gathered}
\label{eq:BigCharDefinitions}
\end{equation}

\begin{remark}
    Equations \eqref{eq:CharSplitEnergyWGamma} and \eqref{eq:SimplifiedForm} show that that even if the boundary data is computed at the correct point ($\gamma = 0$), there are still errors due to errors in the metric terms and Jacobian. It is not enough just to compute the boundary condition values correctly.
\end{remark}
\section{Error Bounds}\label{Sec:ErrorBounds}

Eqns. \eqref{eq:CorrectSolutionBound} and \eqref{eq:ErroneousEnergyBound} show that the correct and erroneous solutions are each bounded, from which it follows that the error is also bounded:
\begin{theorem}
The error $\inorm{\statevec v - \statevec q}_J$ is bounded by data.
\label{thm:BoundedError}
\end{theorem}
\begin{proof}
 By the triangle inequality,
\begin{equation}
\inorm{\statevec v - \statevec q}_J \le \inorm{\statevec v}_J + \inorm{\statevec q}_J.
\end{equation}
The norm of the correct solution is bounded in \eqref{eq:CorrectSolutionBound}. To get the required bound for the erroneous solution,
\begin{equation}
\inorm{\statevec v}_J = \int_\mathcal D J\statevec v^T\statevec v d\xi = \int_\mathcal D \left(\frac{J}{J_e}\right)J_e\statevec v^T\statevec v d\xi = \int_\mathcal D \rho J_e\statevec v^T\statevec v d\xi\le \inorm{\rho}_\infty\inorm{\statevec v}_{J_e}.
\end{equation}
Therefore,
\begin{equation}
\inorm{\statevec e}_J\le \inorm{\rho}_\infty\inorm{\statevec v}_{J_e} + \inorm{\statevec q}_J.
\end{equation}
The ratio of Jacobians, $\rho$, is bounded when the domains are valid, that is, $\rho \le {J_{max}}/{J_{e,min}}$, where $J_{e,min}>0$.
Lastly, $\inorm{\statevec v}_{J_e} $ and $ \inorm{\statevec q}_J$ are each bounded by data through  \eqref{eq:CorrectSolutionBound} and \eqref{eq:ErroneousEnergyBound}.\qed
\end{proof}

Actually, we want to know how the error depends on the geometry errors, so we construct bounds using the error equation, \eqref{eq:CharSplitEnergyWGamma}. For convenience, we define the volume norm of a matrix, $\mmatrix M$, to be
\begin{equation}
\inorm{\mmatrix M}_{2,\infty} = \max_{\mathcal D} \inorm{\mmatrix M}_2,
\end{equation}
where $\inorm{\mmatrix M}_2 = \max_{\statevec w\ne 0}\frac{\statevec w^T\mmatrix M\statevec w}{\statevec w^T\statevec w}$ is the matrix 2-norm. 
Then
\begin{equation}
V1\le \oneHalf\inorm{ \frac{\contraspacevec{\mmatrix E}\cdot\nabla_\xi\varepsilon}{J}}_{2,\infty}\inorm{\statevec e}_{J}^2\equiv c_1 \inorm{\statevec e}_{J}^2
,\quad
V2\le \oneHalf\inorm{ \frac{\contraspacevec{\mmatrix A}\cdot\nabla_\xi\varepsilon}{J}}_{2,\infty}\inorm{\statevec e}_{J}^2 \equiv c_2 \inorm{\statevec e}_{J}^2,
\label{eq:c1c2defs}
\end{equation}
where we have used 
the fact that the matrices are symmetric.
Next, using the Cauchy-Schwarz inequality and the fact that $\contraspacevec{\mmatrix E}$ is divergence-free,
\begin{equation}
V3\le \inorm{\frac{(1+\varepsilon)\contraspacevec{\mmatrix E}\cdot\nabla_\xi\statevec q } {J^2}}_{J}\inorm{\statevec e}_{J} \equiv c_3 \inorm{\statevec e}_{J},\quad V4\le\inorm{\frac{\varepsilon\contraspacevec{\mmatrix A}\cdot\nabla_\xi\statevec q } {J^2}}_{J}\inorm{\statevec e}_{J} \equiv c_4 \inorm{\statevec e}_{J}.
\label{eq:c3c4Defs}
\end{equation}
The coefficients $c_3$ and $c_4$ depend on the correct solution and will serve as source terms for the error.

When we gather the terms in \eqref{eq:SimplifiedForm} and divide by the norm of the error, we get the equation for the main error bound,
\begin{equation}
\ddt \inorm{\statevec e}_{J} +\left(\eta(t) - \gamma\mathcal B_{max}-c_1 - c_2 \right)\inorm{\statevec e}_{J} \le   c_3 + c_4,
\label{eq:MainTimeDerivBound}
\end{equation}
where 
\begin{equation}
\eta(t) = \frac{D}{ \inorm{\statevec e}^{2}_{J}},
\quad \mathcal B =\frac{BG}{ \inorm{\statevec e}^2_{J}},
\end{equation}
and $\mathcal B_{max}$ is the maximum over the time interval.
The difference $\eta(t) - \gamma\mathcal B_{max}$ represents the balance between the advection of error out of the domain and the error introduced into the domain. If negative, error is injected faster than it can be removed through boundary dissipation. The remaining term, $-c_1 - c_2$, depends on $\varepsilon$, the error in the Jacobian. Note that $\inorm{\statevec e}_J = 0$ only if the numerators in $\eta(t)$ and $\mathcal B$ are also zero, since the norm includes the boundary values under the smoothness assumptions.

When $\eta(t) > 0$, Nordstr\"om \cite{doi:10.1137/060654943} argues that for wave propagation problems, $\eta(t)$ is a monotonically growing function over any finite time interval and that the mean value, $\bar\eta$, is bounded from below by a positive constant, that is, $\bar\eta\ge c_0 >0$. Therefore, we can use an integrating factor to get the desired bound.  Assuming that the correct solution, $\statevec q$, is bounded in time so that $c_{3}$, $c_{4}$ are constant independent of time, and defining $\alpha = (c_{0} - \gamma\mathcal B_{max}-c_1 - c_2)$,
\begin{equation}
\begin{split}
\inorm{\statevec e(t)}_{J} &\le e^{-\alpha t}\inorm{\statevec e(0)}_{J} +  (c_{3}+c_{4})\frac{1 - e^{-\alpha t}}{\alpha}.
\end{split}
\label{eq:MainErrorBound}
\end{equation}
Asymptotically, for large time, 
\begin{equation}
\inorm{\statevec e(t)}_{J}\le \frac{1}{\alpha}(c_{3}+c_{4}),
\label{eq:LongTermError}
\end{equation}
when $\alpha \ge 0$, i.e., there is sufficient outflow dissipation relative to the other errors.
Eq. \eqref{eq:LongTermError} relates the long-term error to the geometry errors through the coefficients $c_3$ and $c_4$ defined in \eqref{eq:c3c4Defs}.
\begin{remark}
   The estimates \eqref{eq:MainErrorBound} and \eqref{eq:LongTermError} give qualitative error bounds, due to the unknown quantities $c_0$ and $\mathcal B_{max}$. However, Thm \ref{thm:BoundedError} already tells us that the error is bounded.
\end{remark}

From \eqref{eq:MainTimeDerivBound} and \eqref{eq:MainErrorBound}, we make the following observations:
\begin{enumerate}
\item The coefficient $c_{1}$ is a secondary quantity in that it depends on the product of two (small) error quantities; $c_{2}$ is a primary quantity. Both depend on the ratio of the Jacobians through $\nabla_\xi \varepsilon = \nabla_\xi (J/J_e)$. These terms contribute to growth in the error as a function of time.
\item It is preferred that the exact boundary locations are available at which to evaluate boundary data, i.e. $\gamma = 0$ so $\gamma\mathcal B_{max} = 0$. Nevertheless, $\mathcal B$ varies as $|\Delta\spacevec X|^2$, and so is only a secondary quantity.
\item Even if the boundary data is specified using the exact mapping, $\gamma = 0$, errors represented by $c_3$ and $c_4$ depend on the solution, the error in the Jacobian, $\varepsilon$, and the volume weighted contravariant vectors (through $\contraspacevec{\mmatrix E})$, and contribute to the persistence of the error in time.
\item The dissipation due to outgoing waves ($D$) that comes from using characteristic boundary conditions is necessary to keep the error bounded for long times, as it counteracts the growth due to $c_{1}$, $c_{2}$ and $\mathcal B$, which can be kept small by keeping the erroneous domain sufficiently close to the correct one.
\item \textcolor{black}{Initially, the only error is due to the initial condition. As time evolves, that error is dissipated through the boundaries and the metric errors grow. For long times only the metric errors are significant, cf.~\cite{nordstrom2024uncertain}.}
\item As $\Omega_e\rightarrow \Omega$, $\varepsilon = J/J_e - 1 \rightarrow 0$, $\Delta J\spacevec a^i = J\spacevec a^i_e - J\spacevec a^i\rightarrow 0$ so $c_1, \ldots c_4\rightarrow 0$. For smooth domains, $\nabla_\xi\varepsilon\rightarrow 0$ as $\Omega_e\rightarrow \Omega$, so $\inorm{\statevec e(t)}_{J}\rightarrow 0$. In other words, the mesh error converges to zero as the erroneous domain approaches the correct one.
\end{enumerate}

\subsection{Error Bounds for a Two Dimensional Domain with One Curved Boundary}

\begin{figure}[htbp] 
   \centering
   \includegraphics[width=0.5\textwidth]{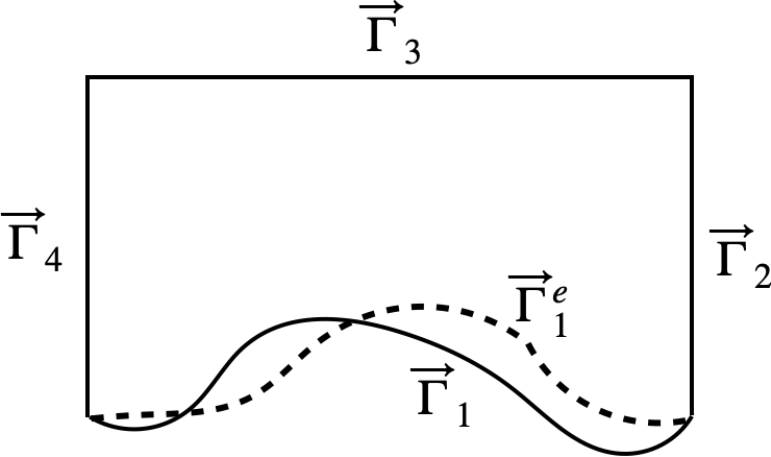} 
   \caption{Correct domain (solid lines) and Erroneous domain (dashed line) boundary curves in two space dimensions used to isolate the boundary errors to one boundary}
   \label{fig:OneCurved.pdf}
\end{figure}

Eq. \eqref{eq:MainTimeDerivBound} is the general result that describes how the error is bounded by data. It is not easy to see, however, how errors in the boundary curves or surfaces directly contribute.
For insight into the direct relationships between the solution error and boundary curve errors, we examine the special case of a simple two-dimensional geometry with one curved boundary, as seen in Fig. \ref{fig:OneCurved.pdf}. A single boundary condition will be applied along the bottom boundary, either specified characteristic data or perfect reflection, \eqref{eq:PerfectReflectionCorrect}.
If more boundary curves have errors, one would simply add their contributions to the mappings.

We construct the two-dimensional mapping $\spacevec X\left(\spacevec \xi\right) = X\left(\xi,\eta\right)\hat x+Y\left(\xi,\eta\right)\hat y + 0\hat z$ for the geometry in Fig. \ref{fig:OneCurved.pdf} as a transfinite interpolation \cite{Gordon&Hall1973a} of the boundaries, as is often done \cite{Kopriva:2009nx}. For the simple geometry of Fig. \ref{fig:OneCurved.pdf}, the correct transfinite mapping is
\begin{equation}
\spacevec X\left(\xi,\eta\right) = \spacevec \Gamma_1(\xi)(1-\eta) + \eta \spacevec \Gamma_3(\xi),
\end{equation}
whereas the erroneous one is
\begin{equation}
\spacevec X^e\left(\xi,\eta\right) = \spacevec \Gamma_1^e(\xi)(1-\eta) + \eta \spacevec \Gamma_3(\xi),
\end{equation}
so that along the lower boundary,
\begin{equation}
\Delta \spacevec X(\xi,0) = \spacevec \Gamma_1^e - \spacevec \Gamma_1  \equiv \Delta \spacevec \Gamma
\label{eq:DeltaX}
\end{equation}
is the position/location error at that boundary. The quantity $\Delta \spacevec \Gamma$ is therefore the error that appears in the boundary terms, \eqref{eq:BoundaryErrorTerm}, for this example. Everywhere else, its influence throughout the domain is 
\begin{equation}
\Delta \spacevec X(\xi,\eta) = (1-\eta)\left(\spacevec \Gamma_1^e(\xi) - \spacevec \Gamma_1(\xi)\right) = (1-\eta)\Delta\spacevec\Gamma.
\end{equation}
Therefore, $\Delta \spacevec X(\xi,\eta)$ is linear in $\Delta \spacevec\Gamma$, and so is the position error for any point in the reference domain.

To get the derivatives of the mapping, and from there the errors in the metric terms and Jacobians, we write the erroneous mapping as the correct mapping plus the error
\begin{equation}
\spacevec X^e = \spacevec X + \Delta\spacevec X.
\end{equation}
The derivatives of the erroneous mapping
\begin{equation}
\begin{gathered}
\spacevec X_\xi^e = \spacevec X_\xi + \Delta\spacevec X_\xi = \spacevec X_\xi +(\spacevec \Gamma_1^{\prime,e} - \spacevec \Gamma'_1)(1-\eta) =  \spacevec X_\xi + (1-\eta)\Delta\spacevec\Gamma' \hfill\\
\spacevec X_\eta^e =  \spacevec X_\eta + \Delta\spacevec X_\eta = \spacevec X_\eta-\Delta \spacevec \Gamma \hfill
\end{gathered}
\label{eq:Metric2D}
\end{equation}
 are also written in terms of the correct ones. 

In two space dimensions the metric terms are
\begin{equation}
\begin{gathered}
J\spacevec a^1 = \spacevec a_2\times \hat z =  \spacevec X_\eta\times \hat z\hfill\\
J\spacevec a^2 = \hat z \times \spacevec a_1 = \hat z \times  \spacevec X_\xi
\end{gathered}
\end{equation}
and the Jacobian is
\begin{equation}
J = \hat z\cdot\left(\spacevec a_1\times\spacevec a_2\right) = \hat z \cdot\left(  \spacevec X_\xi\times\spacevec X_\eta\right).
\end{equation}
So on the erroneous domain,
\begin{equation}
\begin{gathered}
J\spacevec a^{1}_e = J\spacevec a^{1} -  \Delta \spacevec \Gamma\times \hat z\hfill\\
J\spacevec a^{2}_e =J\spacevec a^{2} +  (1-\eta) \hat z \times \Delta \spacevec \Gamma'.\hfill
\end{gathered}
\end{equation}
So we see that the metric term errors depend on $\Delta\spacevec \Gamma$ and its derivative $\Delta\spacevec \Gamma'$,
\begin{equation}
\begin{gathered}
\Delta J\spacevec a^{1} = \hat z\times \Delta \spacevec \Gamma  \hfill\\
\Delta J\spacevec a^{2} = (1-\eta) \hat z \times \Delta \spacevec \Gamma' .\hfill
\end{gathered}
\end{equation}
The vector $J\spacevec a^{2}_e(\xi,0)$ is in the inward normal direction to the erroneous boundary, with the correct outward unit normal being
\begin{equation}
\hat n = \frac{\spacevec n}{|\spacevec n|} = -\frac{ J\spacevec a^{2}}{\left|J\spacevec a^{2}\right|} .
\end{equation}

The error matrices \eqref{eq:ErrorMatrices} are
\begin{equation}
\contravec {\mmatrix E}^i =\Delta J\spacevec a^i\cdot \spacevec {\mmatrix A},
\end{equation}
or, explicitly
\begin{equation}
\begin{gathered}
\contravec {\mmatrix E}^1 = \spacevec {\mmatrix A}\cdot\left( \hat z\times\Delta \spacevec \Gamma\right)  
= \hat z\cdot\left(  \Delta \spacevec \Gamma\times \spacevec {\mmatrix A}\right)  \hfill\\
\contravec {\mmatrix E}^2 = (1-\eta)\spacevec {\mmatrix A}\cdot\left(\hat z \times  \Delta \spacevec \Gamma'\right)
= (1-\eta)\hat z\cdot\left(  \Delta \spacevec \Gamma'\times \spacevec {\mmatrix A}\right).   \hfill
\end{gathered}
\end{equation}
So the error matrices are also linear in $\Delta \spacevec\Gamma$ and $\Delta \spacevec\Gamma'$.
The first, $\contravec {\mmatrix E}^1$, is the error due to location of the boundary. The second, $\contravec {\mmatrix E}^2$, is the error due to the error in the tangent to, or derivative, along the boundary, and also represents the error in the normal there. 

The error in the erroneous Jacobian is  bilinear in $\Delta \spacevec\Gamma$ and $\Delta \spacevec\Gamma'$, for
\begin{equation}
\begin{split}
J_e &=  \hat z \cdot\left(  \left(\spacevec X_\xi +(1-\eta)\Delta\spacevec\Gamma'\right)\times \left(\spacevec X_\eta-\Delta \spacevec \Gamma\right)\right)
\\&
= \hat z\cdot\left\{\left(\spacevec X_\xi\times \spacevec X_\eta \right) - \left( \spacevec X_\xi\times \Delta \spacevec \Gamma \right)+ (1-\eta)\left(\Delta\spacevec\Gamma'\times \spacevec X_\eta \right) - (1-\eta)\left( \Delta\spacevec\Gamma'\times \Delta \spacevec \Gamma\right)  \right\}
\\&
=J +  \hat z\cdot\left\{\left( \Delta \spacevec \Gamma\times \spacevec X_\xi \right)+ (1-\eta)\left(\Delta\spacevec\Gamma'\times \spacevec X_\eta - \Delta\spacevec\Gamma'\times \Delta \spacevec \Gamma\right) \right\}
\\&
\equiv J + \epsilon.
\end{split}
\end{equation}
In terms of the Jacobian error, the ratio of the Jacobians is
\begin{equation}
\rho = \frac{J}{J_e} = \frac{J}{J + \epsilon},
\end{equation}
and so
\begin{equation}
\varepsilon = \rho - 1 =  -\frac{\epsilon}{J + \epsilon}.
\label{eq:Varepsilon2D}
\end{equation}
We see that $\varepsilon$ is nonlinear in the location and derivative errors. However, $\epsilon$ will presumably be small, so we will write
\begin{equation}
\rho  = \frac{1}{1+ \epsilon/J} = 1 - \frac{\epsilon}{J} + O(\epsilon^2).
\end{equation}
To lowest order, $\varepsilon \approx -\epsilon/J$ and is bilinear in $\Delta \spacevec\Gamma$ and $\Delta \spacevec\Gamma'$.

Let us gather the geometric errors and use three unit vectors defined as
\begin{equation}
\hat \alpha =  \frac{\Delta \spacevec \Gamma}{\left|\Delta\spacevec\Gamma\right|},\quad \hat \beta =  \frac{\Delta \spacevec \Gamma'}{\left|\Delta\spacevec\Gamma'\right|}\quad \hat \gamma = \frac{\Delta \spacevec \Gamma''}{\left|\Delta\spacevec\Gamma''\right|}
\end{equation}
to write them as
\begin{equation}
\begin{aligned}
\Delta \spacevec X &=   (1-\eta){\left|\Delta\spacevec\Gamma\right|}\hat\alpha  \\
\contravec {\mmatrix E}^1  &= {\left|\Delta\spacevec\Gamma\right|}\hat z\cdot\left( \hat\alpha\times \spacevec {\mmatrix A}\right) 
 \equiv
{{\left|\Delta\spacevec\Gamma\right|}\hat{\mmatrix A}_1},\hfill\\
\contravec {\mmatrix E}^2 &= (1-\eta)\left|\Delta\spacevec\Gamma'\right| \hat z\cdot\left(\hat\beta\times \spacevec {\mmatrix A} \right)
 \equiv
{(1-\eta)\left|\Delta \spacevec \Gamma'\right|\hat{\mmatrix A}_2}  \\
\epsilon &={ \hat z\cdot\left\{\left( \hat\alpha\times \spacevec X_\xi \right)\left|\Delta\spacevec\Gamma\right|+ 
(1-\eta)\left(\hat\beta\times \spacevec X_\eta \right)\left|\Delta\spacevec\Gamma'\right| 
-(1-\eta)\left|\Delta\spacevec\Gamma'\right|\left|\Delta\spacevec\Gamma\right|\left(\hat\alpha\times \hat\beta\right)  \right\}}
\\&
\equiv {r_1\left|\Delta\spacevec\Gamma\right|+ (1-\eta)r_2\left|\Delta\spacevec\Gamma'\right| + (1-\eta)r_3\left|\Delta\spacevec\Gamma'\right|\left|\Delta\spacevec\Gamma\right|}
\hfill \\
\rho  & =  \frac{J}{J + \epsilon}= 1 - \frac{\epsilon}{J} + O(\epsilon^2) \\
\varepsilon &=-\frac{\epsilon}{J + \epsilon} = -\frac{\epsilon}{J} + O(\epsilon^2) . 
\end{aligned}
\label{eq:ErrorValues}
\end{equation}
Note that the functions used in $\epsilon$,
\begin{equation}
r_1 = \hat z\cdot\left( \hat\alpha \times \spacevec X_\xi\right),\quad
 r_2 = \hat z\cdot \left(\hat\beta\times \spacevec X_\eta \right),\quad 
 r_3 =-\hat z\cdot\left(\hat\alpha\times \hat\beta\right),
 \label{eq:riDefs}
\end{equation}
are bounded by the size of the correct mapping derivatives, $\spacevec X_\xi, \spacevec X_\eta$. The matrices 
\begin{equation}
\hat{\mmatrix A}_1 = \hat z\cdot\left( \hat\alpha\times \spacevec {\mmatrix A}\right) ,
 \quad \hat{\mmatrix A}_2 =  \hat z\cdot\left(\hat\beta\times \spacevec {\mmatrix A} \right)
\end{equation}
are bounded by $\max_{|\hat w|\le 1}\inorm{\hat{\mmatrix A}_i\cdot\hat w}_2$, the norm being the spectral radius. Therefore, all error quantities in \eqref{eq:ErrorValues} are bounded, and to lowest order are proportional to $\left|\Delta\spacevec\Gamma\right|$ and $\left|\Delta\spacevec\Gamma'\right|$.

Thus, we see that as $\Delta\spacevec\Gamma \rightarrow 0$ and $\Delta\spacevec\Gamma' \rightarrow 0$, the geometry errors $\varepsilon\rightarrow 0$ and $\contravec {\mmatrix E}\rightarrow 0$ and so by \eqref{eq:LongTermError}, $\inorm{e}_J\rightarrow 0$.

From the geometric errors, we find how the global solution error depends on errors in the boundary curve for the specified data boundary condition. The perfectly reflecting boundary condition does not introduce additional error (its contribution is zero) but merely reflects it, albeit in the wrong direction, so we need only consider the specified boundary value error. In the process, let us rewrite \eqref{eq:SimplifiedForm} as
\begin{equation}
\begin{split}
\inorm{\statevec e}_{J}\frac{d}{dt}\inorm{\statevec e}_{J}
+ \oneHalf\int_0^{1}  \left. \statevec e^T \mmatrix A^+\statevec e\right|_{\eta=0} d\xi &=-D1 - D2
+\gamma BG
+V1 + V2 + V3 +V4
\\&
\equiv\mathcal R.
\end{split}
\label{eq:CharSplitEnergy2D2}
\end{equation}
The non-negative quantities $D1 = \oneHalf\int_0^{1}  \varepsilon\left. \statevec e^T \mmatrix A^+\statevec e\right|_{\eta=0} d\xi$, $D2=\oneHalf\int_{0}^{1}  \left. (1+\varepsilon)\statevec e^T \mmatrix E^+\statevec e\right|_{\eta=0} d\xi$, and $BG$ are the dissipation and boundary errors whose dependence on the geometry errors is derived in \eqref{eq:D1}, \eqref{eq:BoundarydissipationFnofDelta}, \eqref{eq:DG}. In \eqref{eq:CharSplitEnergy2D2} have used the fact that the boundary errors are nonzero only along the $\spacevec \Gamma_{1}= \spacevec\Gamma$ boundary.

The functional, $\mathcal R$, contains all the error terms controllable by the accuracy of $\spacevec \Gamma_1$.
In \ref{AppB}, we show that
\begin{equation}
\begin{split}
|\mathcal R(\Delta\spacevec\Gamma, \Delta\spacevec\Gamma',\Delta\spacevec\Gamma'')| 
&\le  \left|\Delta\spacevec\Gamma\right|_{max}\left\{\oneHalf \inorm{\frac{\mmatrix M_1}{J}}_{2,\infty } \inorm{ \statevec e}_J 
+\inorm{\frac{\hat{\mmatrix A}_1\statevec q_\xi}{J^2}}_J 
+\inorm{\frac{r_1\left(\contraspacevec{\mmatrix A}\cdot\nabla_\xi\statevec q\right)}{J^3}}_J
\right\}\inorm{ \statevec e}_J
 \\&+ | \Delta \spacevec \Gamma'|_{max}\left\{\oneHalf\inorm{\frac{\mmatrix M_2}{J}}_{2,\infty }\inorm{ \statevec e}_J  
  +\inorm{\frac{(1-\eta)\hat{\mmatrix A}_2\statevec q_\eta}{J^2}}_J +
\inorm{\frac{(1-\eta)r_2\left(\contraspacevec{\mmatrix A}\cdot\nabla_\xi\statevec q\right)}{J^3}}_J
   \right\} \inorm{ \statevec e}_J \\&
 + | \Delta \spacevec \Gamma''|_{max}\left\{\oneHalf\inorm{\frac{\mmatrix M_3}{J}}_{2,\infty }\inorm{ \statevec e}_J  
   \right\} \inorm{ \statevec e}_J + HOTs,
\end{split}
\label{eqeq:BigRBoundRepeat}
\end{equation}
where the matrices $\mmatrix M_1, \mmatrix M_2,\mmatrix M_3$ are defined in \eqref{eq:Mmatrices}.
The first term in each brace in \eqref{eqeq:BigRBoundRepeat} contributes to the exponential behavior of the error in time through the coefficient $\alpha$ in \eqref{eq:MainErrorBound}. The remaining terms depend on the correct solution and act as source terms, affecting the size of the long term error through $c_3$ and $c_4$. The $\mathcal B_{max}$ term in \eqref{eq:MainTimeDerivBound} is second and higher order.

The long term error depends on the size of $c_3$ and $c_4$, see \eqref{eq:LongTermError}. Using the definitions for $c_3,c_4$ in \eqref{eq:c3c4Defs}, and the bounds from \eqref{eq:V3AppBBound} and \eqref{eq:V4AppBBound} in \ref{AppB},
\begin{equation}
\begin{split}
c_3 + c_4 &
= \left|\Delta\spacevec\Gamma\right|_{max}
\left\{\inorm{\frac{\hat{\mmatrix A}_1\statevec q_\xi}{J^2}}_J
+ \inorm{\frac{r_1\left(\contraspacevec{\mmatrix A}\cdot\nabla_\xi\statevec q\right)}{J^3}}_J\right\} 
\\&
+  \left|\Delta\spacevec\Gamma'\right|_{max}
\left\{\inorm{\frac{(1-\eta)\hat{\mmatrix A}_2\statevec q_\eta}{J^2}}_J + \inorm{\frac{(1-\eta)r_2\left(\contraspacevec{\mmatrix A}\cdot\nabla_\xi\statevec q\right)}{J^3}}_J \right\} + HOTs.
\end{split}
\label{eq:c3+c4Coefs}
\end{equation}
From \eqref{eq:c3+c4Coefs} and \eqref{eq:LongTermError}, we conclude that for long times to the lowest order, the global error is proportional to the magnitude of the location error $\Delta \spacevec \Gamma$ and to the magnitude of the derivative error $\Delta \spacevec \Gamma'$. The coefficient of $\left|\Delta\spacevec\Gamma'\right|_{max}$ depends on $(1-\eta)$, which vanishes at the upper boundary. All things being equal, that factor makes those coefficients smaller than the coefficients of $\left|\Delta\spacevec\Gamma\right|_{max}$ for the simple geometry in Fig. \ref{fig:OneCurved.pdf}. Finally, then, from \eqref{eq:LongTermError}, we see that the error converges linearly with $\Delta \spacevec \Gamma$ {\it and} $\Delta \spacevec \Gamma'$ as $\Omega_e\rightarrow \Omega$.
\begin{remark}
    The dependence, $1-\eta$, in the influence of the derivative error on the solution error is a result of the transfinite interpolation with linear blending for the mapping, and reflects the fact that the bottom boundary does not affect top one.
\end{remark}
%
%
\section{Examples}
We present two examples to demonstrate the predicted linear variation \eqref{eq:c3+c4Coefs} of the error with respect to boundary location and derivative errors. The first, in one space dimension, can be studied from the analytic solution, but contains only boundary location error. The second, in two space dimensions, has the single boundary curve error as in Fig. \ref{fig:OneCurved.pdf}.

\subsection{The Scalar Advection Equation in one Space Dimension}

The first example solves the IBVP for the constant coefficient scalar advection equation
\begin{equation}
\begin{gathered}
q_t + a q_x = 0,\quad x\in\Omega = [0,1], a>0\hfill\\
q(0,t) = -\sin(2\pi at)\hfill\\
q(x,0) = \sin(2\pi x).\hfill
\end{gathered}
\label{eq:1DScalarProblem}
\end{equation}
In this case, $\Omega = \mathcal D$, $J=1$, and $\xi = x$. 
The correct solution is $q(\xi,t) = \sin(2\pi(\xi-at))$.

Now let $\Omega_e = [\delta,1]$, $\delta<1$, be the erroneous domain. Transforming to the unit reference domain,
\begin{equation}
X_e = \delta + \xi(1-\delta), \quad J_e = \frac{\partial x}{\partial \xi} = 1-\delta,
\end{equation}
so that from \eqref{eq:DeltaX}, $\Delta X(\xi) = (1-\xi)\delta$, and from \eqref{eq:Metric2D}, $J\spacevec a^1 = 1\hat x$.
Then the erroneous domain problem becomes
\begin{equation}\begin{gathered}
v_t + \frac{a}{(1-\delta)} v_\xi = 0 \hfill\\
v(0,t) = g_\gamma(\xi,t)\hfill\\
v(x,0) = \sin\left(2\pi(\delta + \xi(1-\delta))\right).\hfill
\end{gathered}
\end{equation}
We consider two choices for the boundary value $g_\gamma(\xi,t)$, namely
\begin{equation}
        g_\gamma = \sin(2\pi(\gamma\delta - at)),
\end{equation}
with $\gamma =0$ and $\gamma=1$.
As used above in \eqref{eq:GammaDef}, $g_0$ is exact at the left boundary, and no error is introduced there. For $g_1$, the boundary value is evaluated at the erroneous location, and hence introduces error there.
The analytic solution for the erroneous domain is then
\begin{equation}
v(\xi,t) =
\left\{
\begin{gathered}
 \sin\left(2\pi\left(\delta+\xi(1-\delta) - \frac{a}{1-\delta}t\right)\right)\quad \xi - \frac{a}{1-\delta}t>0 \\
 -\sin\left(2\pi a\left(t - \frac{\gamma\delta + \xi(1-\delta)}{a}\right)\right)\quad \xi - \frac{a}{1-\delta}t\le0\hfill
 \end{gathered}
\right. .
\end{equation}

In Sec. \ref{Sec:ErrorBounds} we showed what quantities affect the error and its bounds. For the one dimensional scalar problem here,
\begin{equation}
\varepsilon = \frac{1}{1-\delta} - 1 = \frac{\delta}{1-\delta} = \delta\rho,
\end{equation} 
so that $\nabla_\xi\varepsilon = 0$. Next, $\contraspacevec{\mmatrix E} = 0$ since $J\spacevec a^1 = J\spacevec a^1_e = 1\hat x$.  Then $c_1 = c_2 = 0$ and $c_3=0$ in \eqref{eq:MainTimeDerivBound}. But $c_4 = 2\pi|\delta/(1-\delta)|$, since the correct solution is sinusoidal. When $\gamma = 0$, $e_g=0$, and the inflow boundary term, $\mathcal B$, vanishes. Then with the dissipation, $D$, included at the right boundary,
\begin{equation}
\ddt\inorm{ e(t)}_{J} +\left[ \frac{ ae^2(1,t)}{\inorm{ e(t)}^2_{J}}\right] \inorm{ e(t)}_{J}\le  \frac{4\pi|\delta|}{1-\delta} = 4\pi|\delta|\rho,
\label{eq:Scalar1DErrorBound}
\end{equation}
which, since $\rho$ is a constant, we write as
\begin{equation}
\ddt\inorm{ e(t)}_{J} +\rho\eta(t) \inorm{ e(t)}_{J}\le  4\pi|\delta|\rho.
\label{eq:Scalar1DErrorBound2}
\end{equation}
As argued above, the mean value of $\eta(t)$ over any finite time interval, $\bar\eta$, is bounded from below by a positive constant, i.e., $\bar\eta\ge c_0 >0$ and therefore we can use an integrating factor to the bound
\begin{equation}
\inorm{ e(t)}_{J} \le e^{-\rho c_0t}\inorm{ e(0)}_{J} + \frac{1 - e^{-\rho c_0t}}{\rho c_0} 4\pi|\delta|\rho =e^{-\rho c_0t}\inorm{ e(0)}_{J} + \frac{1 - e^{-\rho c_0t}}{ c_0} 4\pi|\delta|.
\label{eq:1DErrorBoundInTime}
\end{equation}
As seen in Sec. \ref{Sec:ErrorBounds} for the general case, \eqref{eq:1DErrorBoundInTime} says that the error is bounded in time, for large time it approaches $4\pi|\delta|/c_0$, and it varies linearly with the boundary location error, $\delta$. We also argued in Sec. \ref{Sec:ErrorBounds} that the boundary value evaluation error introduced when applying the boundary condition at the wrong point in space ($\gamma = 1$) should not change this result significantly, since that term is a secondary quantity. In other words, the dominant error should be due to the volume errors.

\begin{figure}[htbp] 
   \centering
   \includegraphics[width=0.7\textwidth]{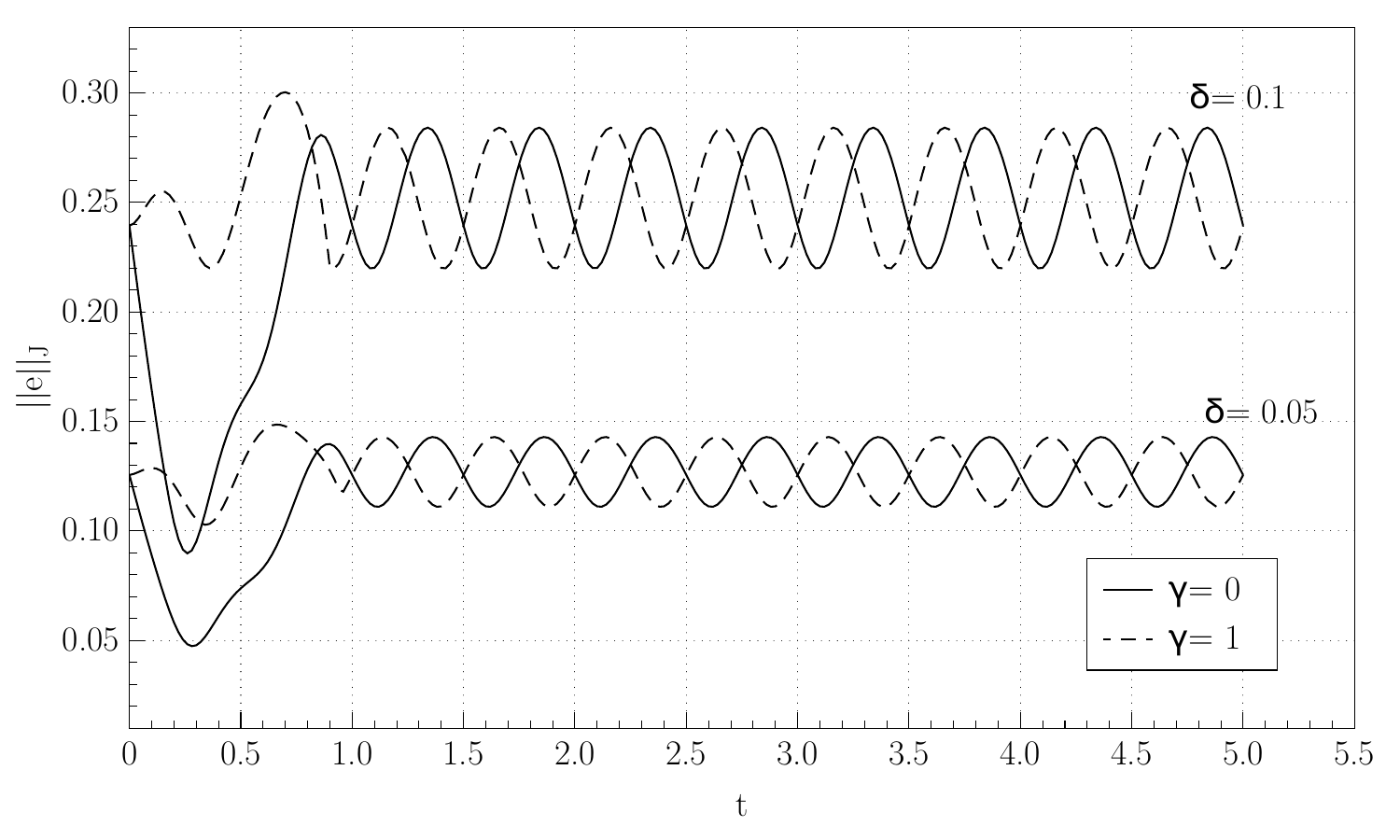} 
   \caption{\textcolor{black}{Time history of the error for two values of the left boundary location error, $\delta$, and for two values of $\gamma$ for the one dimensional scalar problem, \eqref{eq:1DScalarProblem}. The time histories show that the errors are bounded and linearly dependent on $\delta$. For early times the errors differ between $\gamma = 1$ (evaluating boundary conditions at the erroneous boundary location) and $\gamma = 0$ (evaluating at the correct boundary location). For long times, the error bounds are the same, as described in \eqref{eq:1DErrorBoundInTime}. This also illustrates observation 5 from Sec.~\ref{Sec:ErrorBounds}}}
   \label{fig:TimeHistoryError.pdf}
\end{figure}

\begin{figure}[htbp] 
   \centering
   \includegraphics[width=0.7\textwidth]{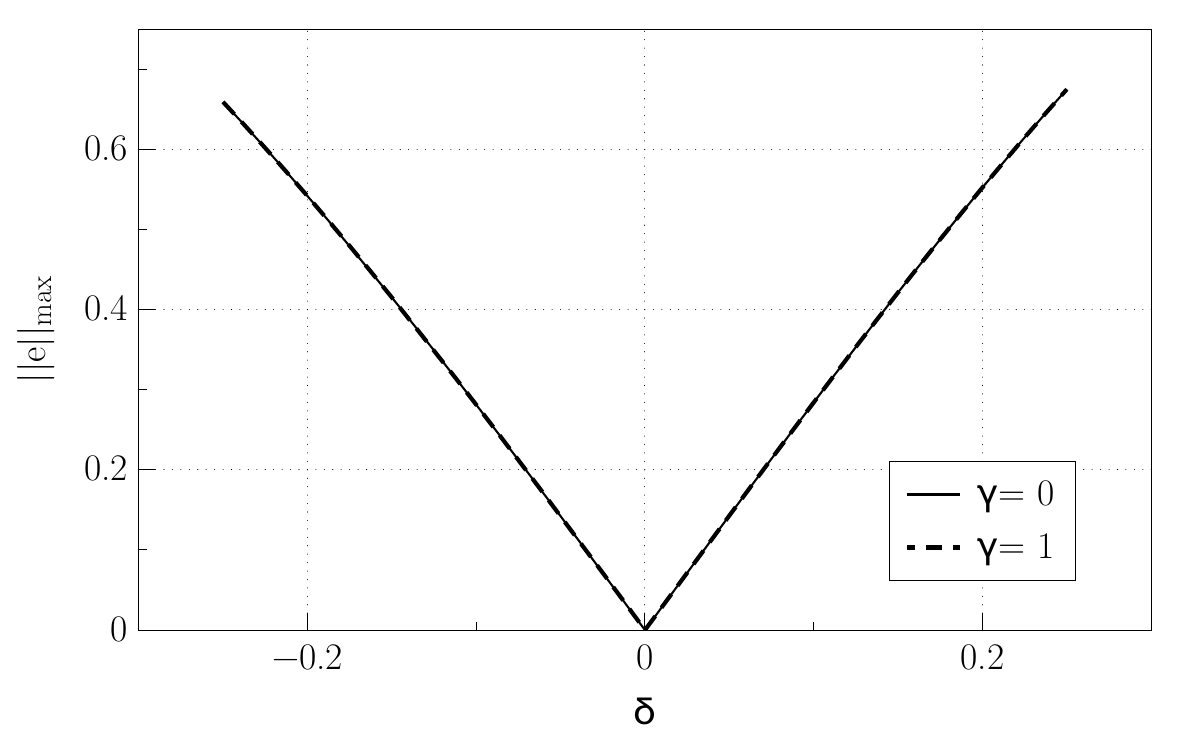} 
   \caption{\textcolor{black}{Variation of the maximum errors as seen in Fig. \ref{fig:TimeHistoryError.pdf} of the large time solution of the one dimensional scalar problem, \eqref{eq:1DScalarProblem} plotted as a function of $\delta$, demonstrating the linear dependence of the error on $\delta$ in \eqref{eq:1DErrorBoundInTime}. The results show the long term minimal effect of evaluating the boundary conditions at the erroneous location ($\gamma = 1$)}}
   \label{fig:EAsFunOfDeltaG}
\end{figure}

When we measure the norm of the error, it is bounded, as predicted, for both $\gamma = 0$ and $\gamma = 1$. Fig. \ref{fig:TimeHistoryError.pdf} shows the time history of the norm of the error for two values of $\delta$ and for both values of $\gamma$. As expected, the error is larger for the larger boundary location error, 0.283 for $\delta = 0.1$ and 0.143 for $\delta = 0.05$. The ratio of the two, 1.98, confirms the linear variation with $\delta$. Further confirmation of the linear behavior of the maximum error for small $|\delta|$ is shown in Fig. \ref{fig:EAsFunOfDeltaG}, which plots the peak value of the error for large time seen in Fig.  \ref{fig:TimeHistoryError.pdf} as a function of $\delta$ for both positive and negative values of $\delta$. Finally, Figs. \ref{fig:TimeHistoryError.pdf} and \ref{fig:EAsFunOfDeltaG} also confirm that the boundary value error, $\mathcal B$, does not significantly impact the maximum of the $L^2$ norm of the error for long times, since the maximum errors are indistinguishable between $\gamma = 0$ and $\gamma = 1$. 


\subsection{The Scalar Advection Equation in Two Space Dimensions}

We illustrate the errors in two space dimensions for the correct domain $\Omega = [0,1]^2$ and erroneous domains with
\begin{equation}
\spacevec X\left(\xi,\eta\right) = \spacevec \Gamma_1(\xi)(1-\eta) + \eta \left[\xi\hat x + 1\hat y \right].
\end{equation} 

\begin{figure}[htbp] 
   \centering
   \includegraphics[width=0.5\textwidth]{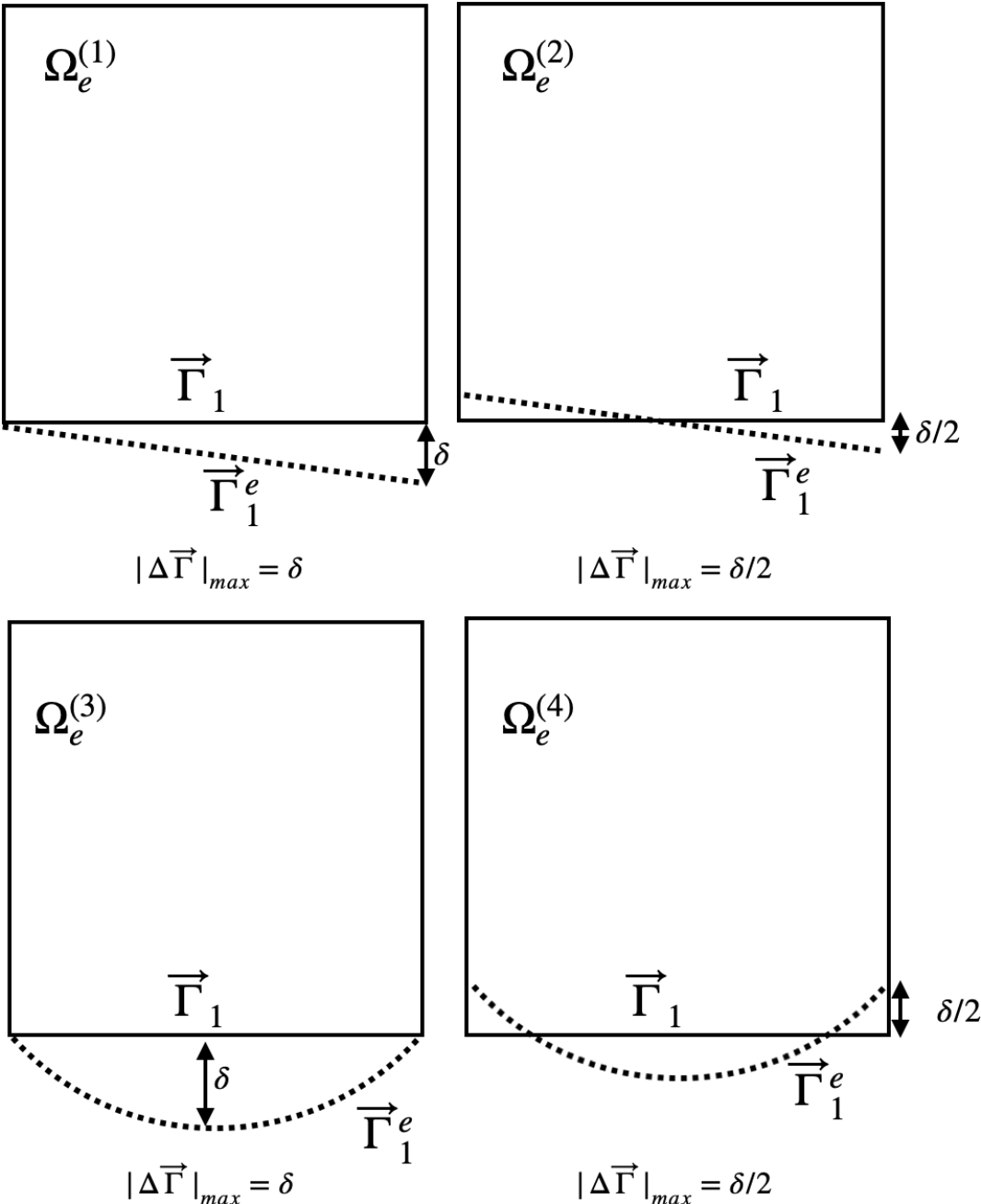} 
   \caption{\textcolor{black}{Correct and erroneous domains with the linear and quadratic perturbations of \eqref{eq:Four2DDomainCurves} of the correct unit square domain used to test \eqref{eq:LongTermError} in two space dimensions.  In the left column, the erroneous curve matches at least one end point. The right column shows a minimax approximation of the bottom boundary}}
   \label{fig:ErroneousAndCorrectDomains}
\end{figure}

We choose four erroneous domains, two with a linear lower boundary and two with a quadratic one, see Fig. \ref{fig:ErroneousAndCorrectDomains}. For two domains, one or two points match the endpoints of the correct bottom curve. The other two are the minimax polynomial approximations of degree one or two. The bottom curves for the four cases are
\begin{equation}
\begin{gathered}
\Omega_e^{(1)} : \Gamma^e_1(\xi) = \xi\hat x - \delta\xi \hat y\hfill\\
\Omega_e^{(2)} : \Gamma^e_1(\xi) = \xi\hat x  - \delta\left(\xi - \oneHalf\right) \hat y\hfill\\
\Omega_e^{(3)} : \Gamma^e_1(\xi) = \xi\hat x + 4\delta\xi(\xi-1)\hat y\hfill\\
\Omega_e^{(4)} : \Gamma^e_1(\xi) = \xi\hat x + \left(4\delta\xi(\xi-1) + \oneHalf\delta\right) \hat y.\hfill\\
\end{gathered}
\label{eq:Four2DDomainCurves}
\end{equation}
The parameter $\delta$ controls the deviation of the erroneous domain from the correct one, and also the error in the derivative. 

The IBVP solved is
\begin{equation}
\begin{gathered}
q_t + \spacevec a\cdot \nabla_x q = 0,\quad \spacevec a = a_1\hat x + a_2\hat y > 0, \quad x\in\Omega, t>0\hfill\\
q(0,y,t) = f(a_2y- |a|^2 t) \hfill\\
q(x,y,t) = f(\spacevec a\cdot\spacevec \Gamma_1 - |a|^2t), \quad \spacevec x\in\spacevec \Gamma_1\hfill\\
q(x,y,0) = f(a_1 x+ a_2 y),\hfill
\end{gathered}
\end{equation}
which has the analytic solution
\begin{equation}
q(x,y,t) = f(\spacevec a\cdot\spacevec x - |a|^2 t).
\label{eq:fourDomainIC}
\end{equation}
For the examples here we choose $f(\phi) = \sin(\omega\pi\phi)$ with $\omega = 4$, and $\spacevec a = (\sqrt{3}/2,1/2)$ so that $|a| = 1$. 

\begin{remark}
The choice of $f$ and its parameters affect only the size of the coefficients $c_3$ and $c_4$ through the bounds on the solution gradient. It does not affect the dependence of the error on $\Delta\spacevec \Gamma$ and $\Delta\spacevec \Gamma'$, see \eqref{eq:c3+c4Coefs}.
\end{remark}
Since this problem becomes variable coefficient in the erroneous domain, we compute its solution numerically with a single domain discontinuous Galerkin Spectral Element approximation, with a third order Runge-Kutta approximation in time. See \cite[Sec. 7.4]{Kopriva:2009nx} for complete details. We compute the solutions so that the approximation errors are small compared to the erroneous domain errors. We choose a polynomial approximation order of $N=18$ and time step $\Delta t=1\times 10^{-4}$ so that the approximation errors were $3.2\times 10^{-8}$ for the correct domain and less than $6\times 10^{-6}$ for any of the erroneous domains at the chosen final time, $t = 1.5$. Finally, since the boundaries are (at most) polynomials of degree two, the boundary mappings and metric terms are represented exactly by the approximations. The examples can therefore be interpreted as measuring $e_{mesh}$, and as a case where $e_{mesh} > e_{solver}$.

Eqs. \eqref{eq:LongTermError} and \eqref{eq:c3+c4Coefs} imply that the domain error should vary linearly with $|\Delta\Gamma |_{max}$ for small $\delta$. For the linear domain, $\Omega^{(1)}_e$ (see Fig. \ref{fig:ErroneousAndCorrectDomains}),
\begin{equation}
|\Delta\spacevec \Gamma|_{max} = |\delta|,\quad |\Delta\spacevec \Gamma'|_{max} = \oneHalf|\delta |,\quad |\Delta\spacevec \Gamma''|_{max} = 0
\end{equation}
and for the quadratic,
\begin{equation}
|\Delta\spacevec \Gamma|_{max} = |\delta|,\quad |\Delta\spacevec \Gamma'|_{max} = |\delta |,\quad |\Delta\spacevec \Gamma''|_{max} = 2|\delta|.
\end{equation}
For the minimax approximations the $|\Delta\spacevec \Gamma|_{max}$ quantities are halved, but the others remain the same. Thus, as long at the coefficients for the perturbation quantities in \eqref{eq:BigRBound} remain constant, the errors should be linear in $\delta$. Furthermore, for each of these erroneous domains, $r_2 = 0$ (see \eqref{eq:riDefs}), making the $|\Delta\spacevec \Gamma|_{max}$ error term dominant in \eqref{eqeq:BigRBoundRepeat}, so the error for the minimax approximations should be roughly half of the others.

\begin{figure}[htbp] 
   \centering
   \includegraphics[width=0.495\textwidth]{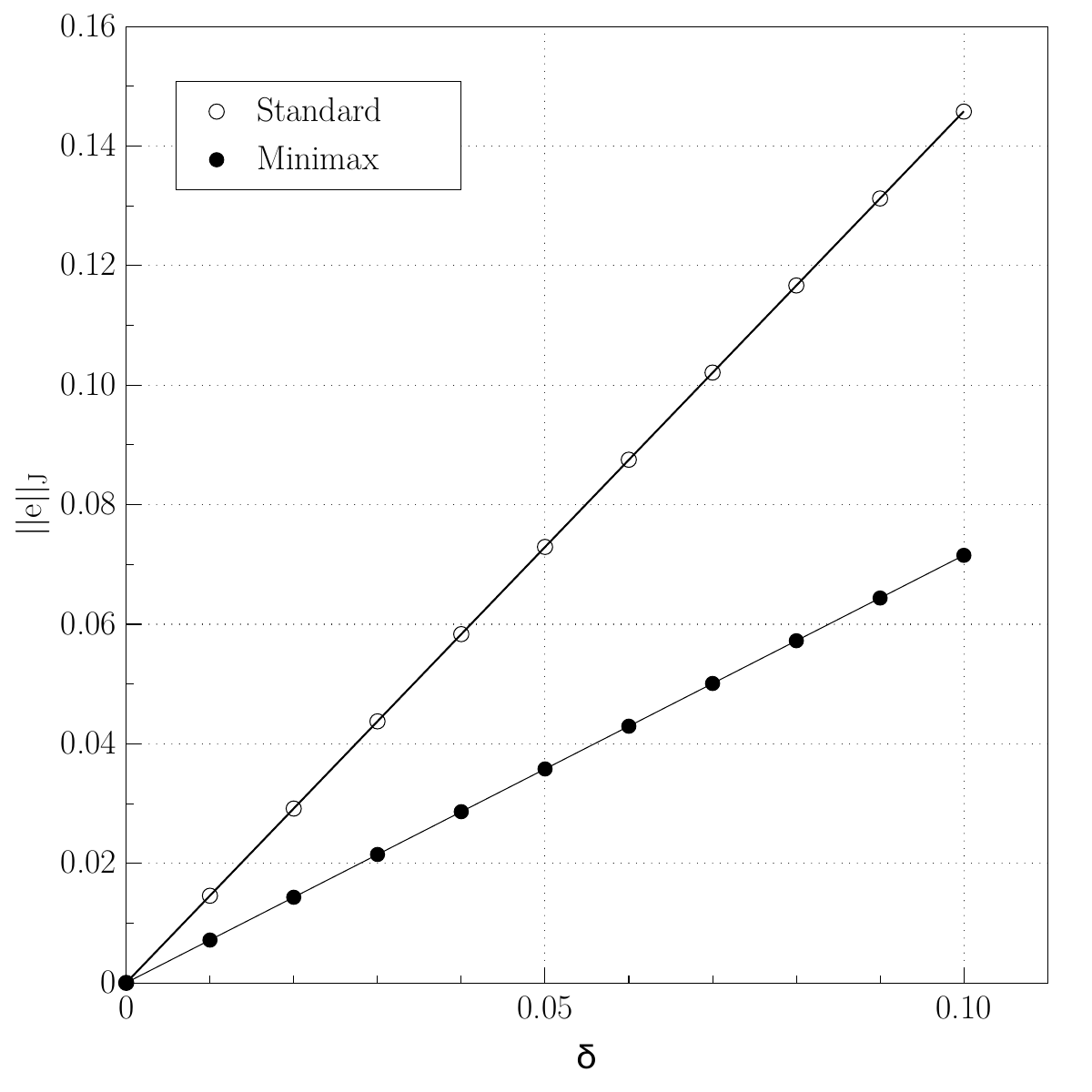} 
   \includegraphics[width=0.495\textwidth]{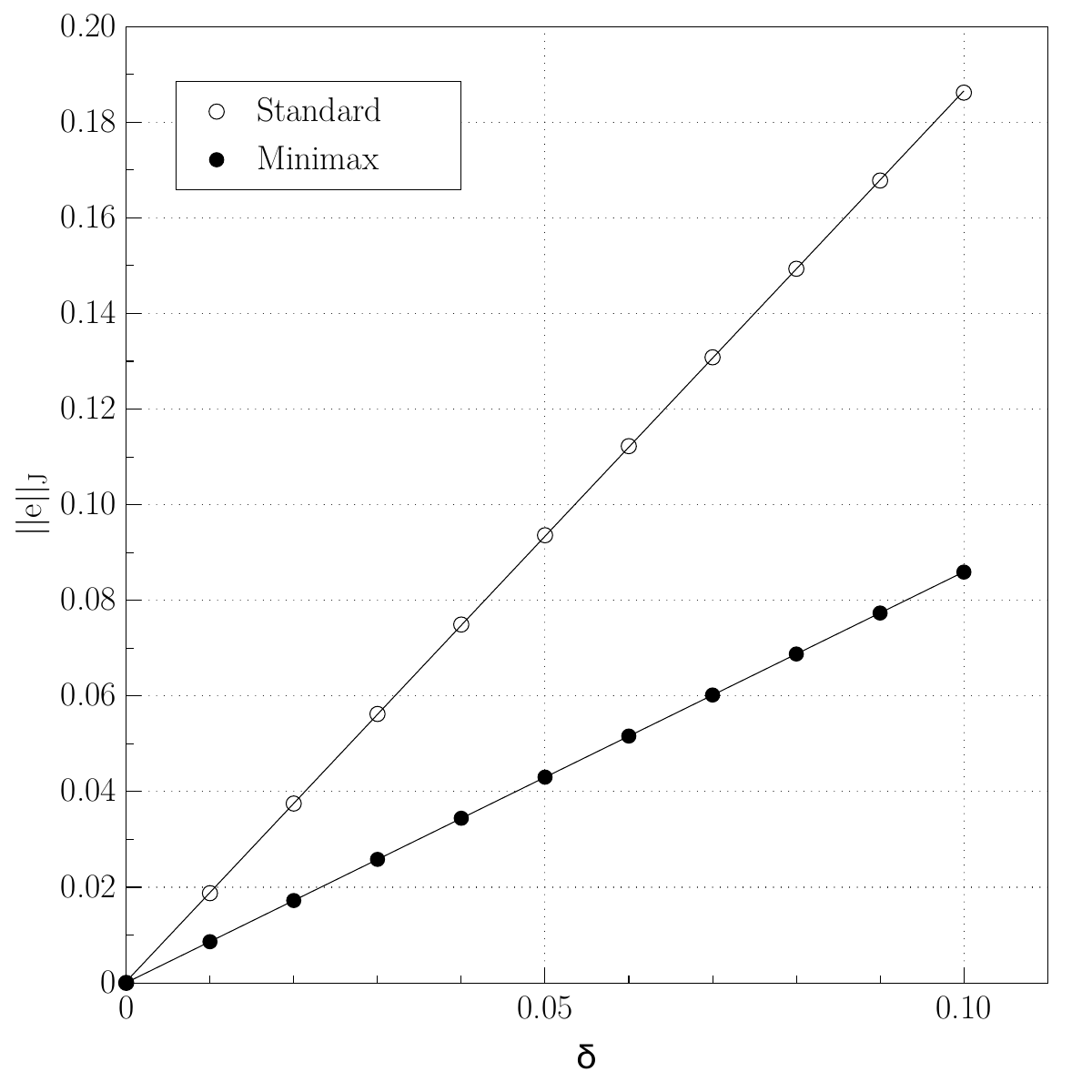} 
   \caption{\textcolor{black}{Linear growth of the erroneous domain error as a function of perturbation parameter $\delta$ for the standard ($\Omega^{(1)}_e$) and minimax ($\Omega^{(2)}_e$) linear erroneous domains (left) and for the ($\Omega^{(3)}_e$) and minimax ($\Omega^{(4)}_e$) quadratic erroneous domains (right) as expected from \eqref{eq:LongTermError} and \eqref{eq:c3+c4Coefs}. Circles are the measured values and the lines are the linear least squares approximations through those points}}
   \label{fig:2DDeltaLinearError}
\end{figure}

In Fig. \ref{fig:2DDeltaLinearError} we plot the error as a function of $\delta$ for the linear erroneous domains and see that it grows linearly, as predicted. The ratio of the slopes is 2.04, which indicates that the error is dominated by the $|\Delta\spacevec \Gamma|_{max}$ term.  We can evaluate the bounds for the coefficients in \eqref{eq:c3+c4Coefs} \textcolor{black}{using the parameters for the problem \eqref{eq:Four2DDomainCurves} - \eqref{eq:fourDomainIC}}, and find that $R\le 25.4|\Delta\Gamma|_{max} + 4.52|\Delta\Gamma'|_{max} + HOTs$ for $\Omega^{(1)}_e$.
In Fig. \ref{fig:2DDeltaLinearError} we also plot the error as a function of $\delta$ for the quadratic erroneous domains and see that it, too, grows linearly, as expected. The ratio of the slopes is 2.18, which also indicates that the error is dominated by the $|\Delta\spacevec \Gamma|_{max}$ term. When we compute bounds for the coefficients in \eqref{eq:c3+c4Coefs}, we find that $R\le 27.2|\Delta\Gamma|_{max} + 3.9|\Delta\Gamma'|_{max} + HOTs$ for $\Omega^{(3)}_e$.

\begin{figure}[htbp] 
   \centering
   \includegraphics[width=0.495\textwidth]{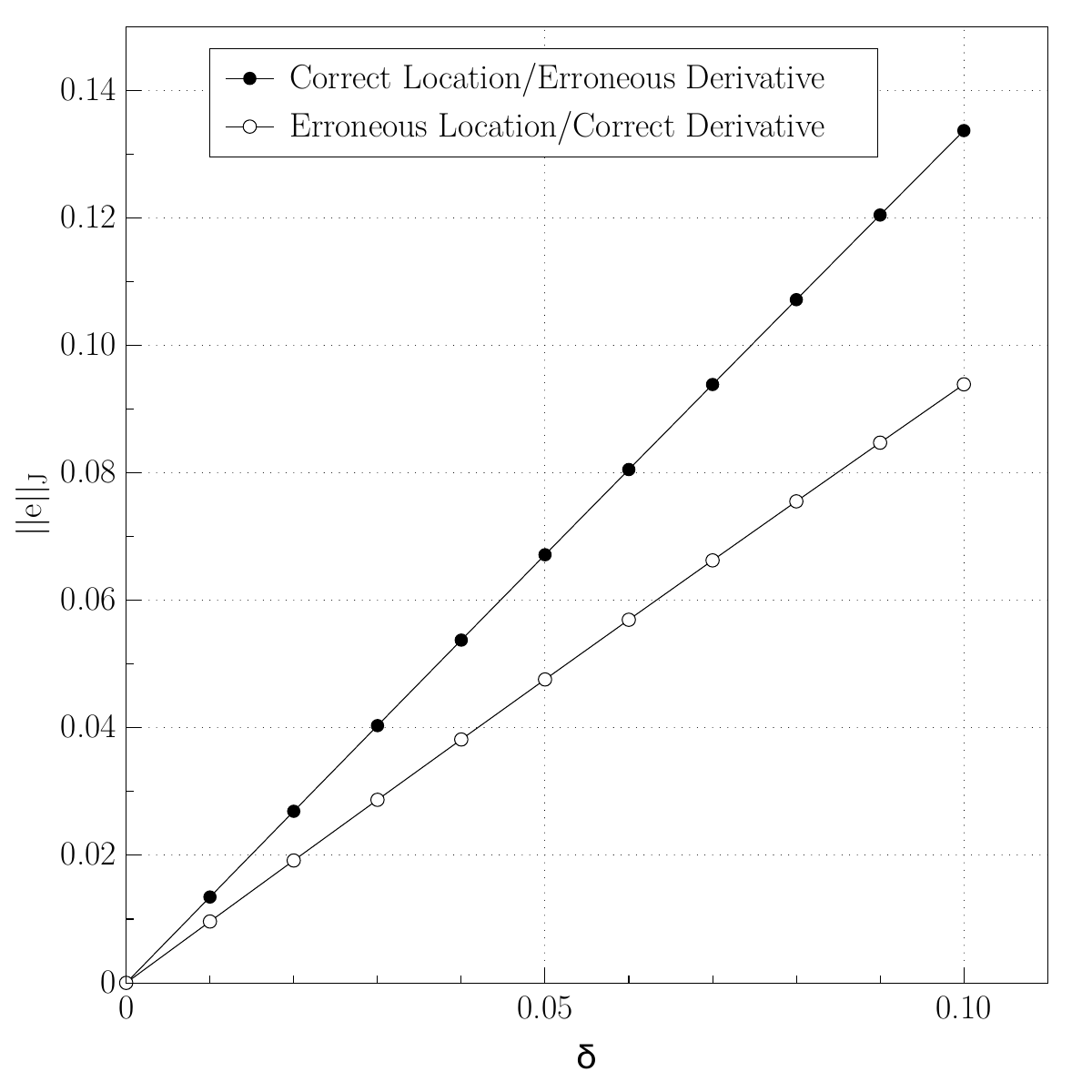} 
    \includegraphics[width=0.495\textwidth]{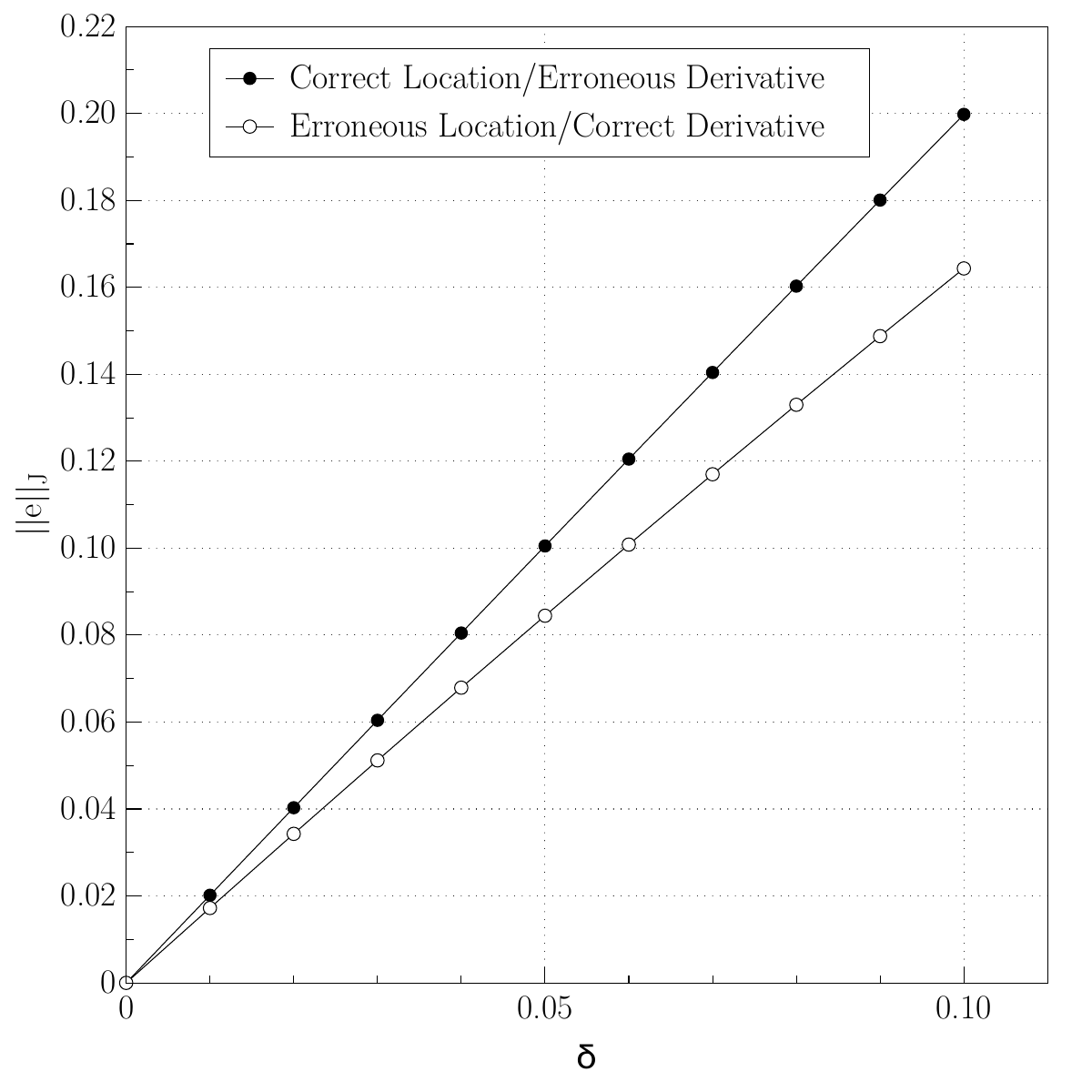} 
  \caption{\textcolor{black}{Isolated errors as a function of $\delta$ for $\Omega^{(1)}_e$ (left) and $\Omega^{(3)}_e$ (right) showing the linear variation of the error for each term in \eqref{eq:c3+c4Coefs}}}
   \label{fig:ECPlots}
\end{figure}

\begin{figure}[htbp] 
   \centering
   \includegraphics[width=0.495\textwidth]{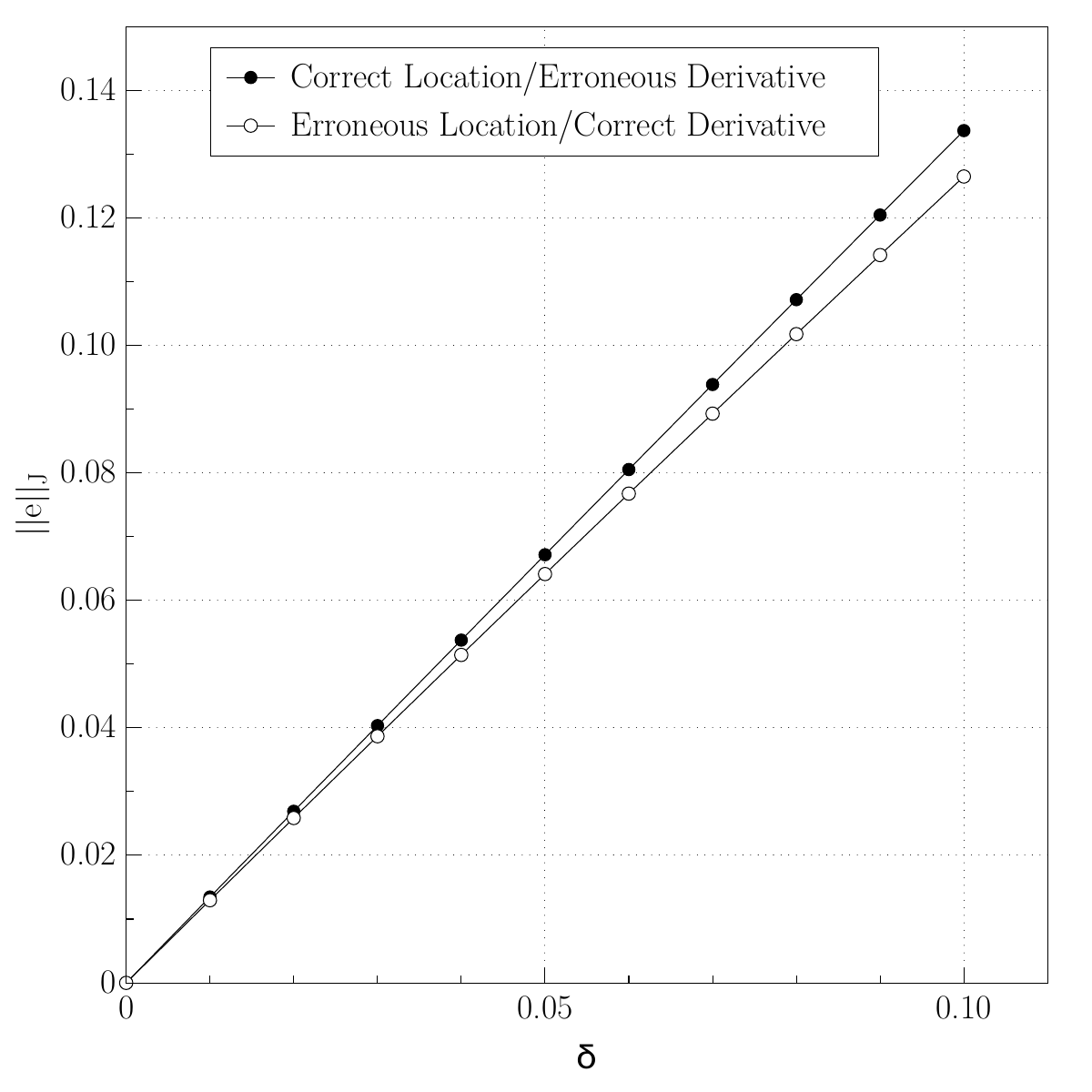} 
    \includegraphics[width=0.495\textwidth]{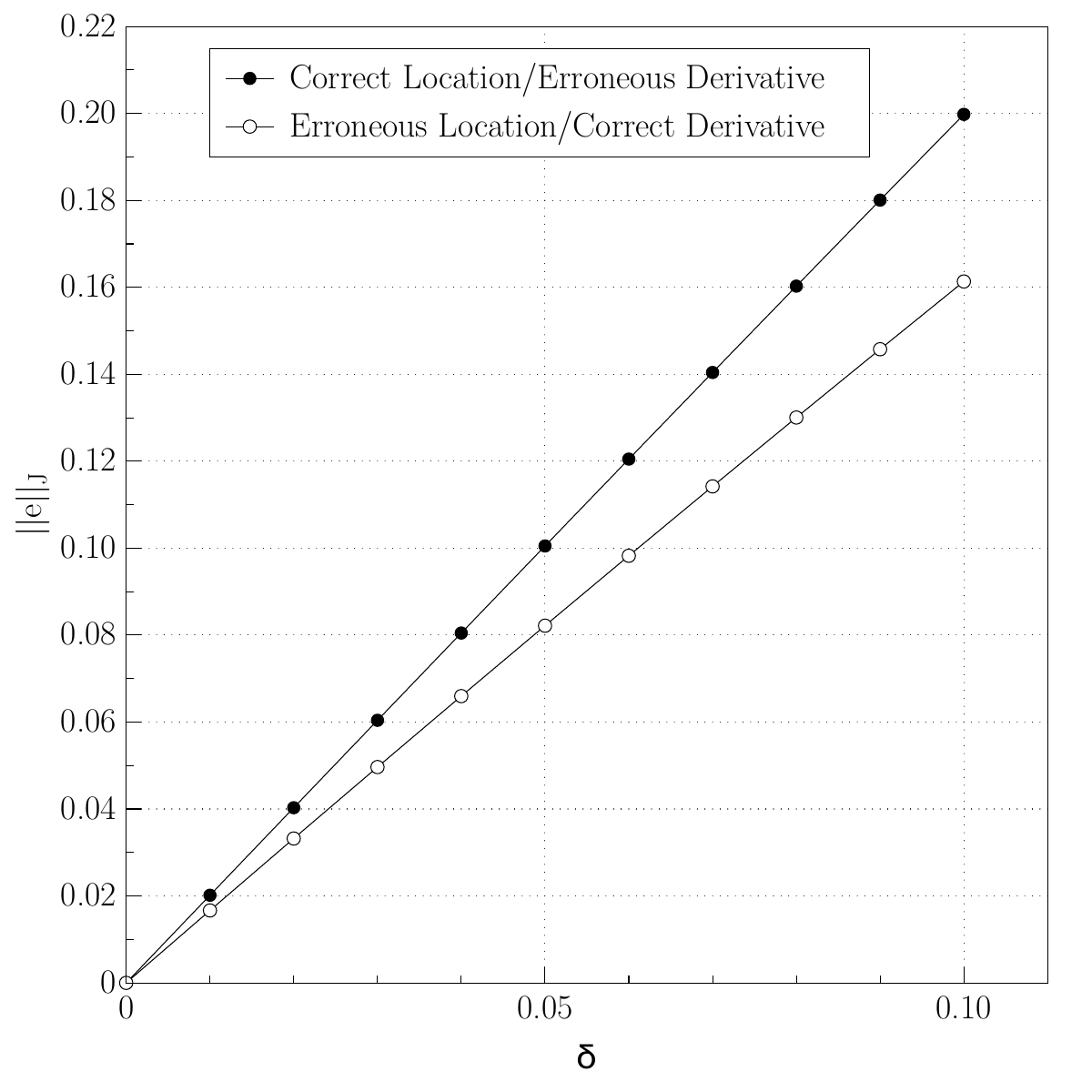} 
  \caption{\textcolor{black}{Isolated errors as a function of $\delta$ for the minimax approximations, $\Omega^{(2)}_e$ (left) and $\Omega^{(4)}_e$ (right) showing the linear variation of the error for each term in \eqref{eq:c3+c4Coefs}} }
   \label{fig:ECPlots24}
\end{figure}

We can isolate the effects of the errors in the location and the derivative in \eqref{eq:c3+c4Coefs} by using the correct values for one and the erroneous values for the other when constructing the geometry and metric terms. For instance, ``Correct Location/Erroneous Derivative" corresponds to $\Delta\spacevec\Gamma = 0$ and $\Delta\spacevec\Gamma' = -\delta$ for $\Omega^{(1)}_e$. In each case, the error should vary linearly in the parameter $\delta$. Fig. \ref{fig:ECPlots} shows that the error does vary linearly for $\Omega_e^{(1)}$ and $\Omega_e^{(3)}$. Fig. \ref{fig:ECPlots24} shows the same for the minimax approximations $\Omega_e^{(2)}$ and $\Omega_e^{(4)}$. We note that for the isolated errors, the errors in the derivative have a significant impact on the solution error, with the slopes differing by up to $42\%$ for the examples presented.


\begin{figure}[htbp] 
  \centering
   \includegraphics[width=0.495\textwidth]{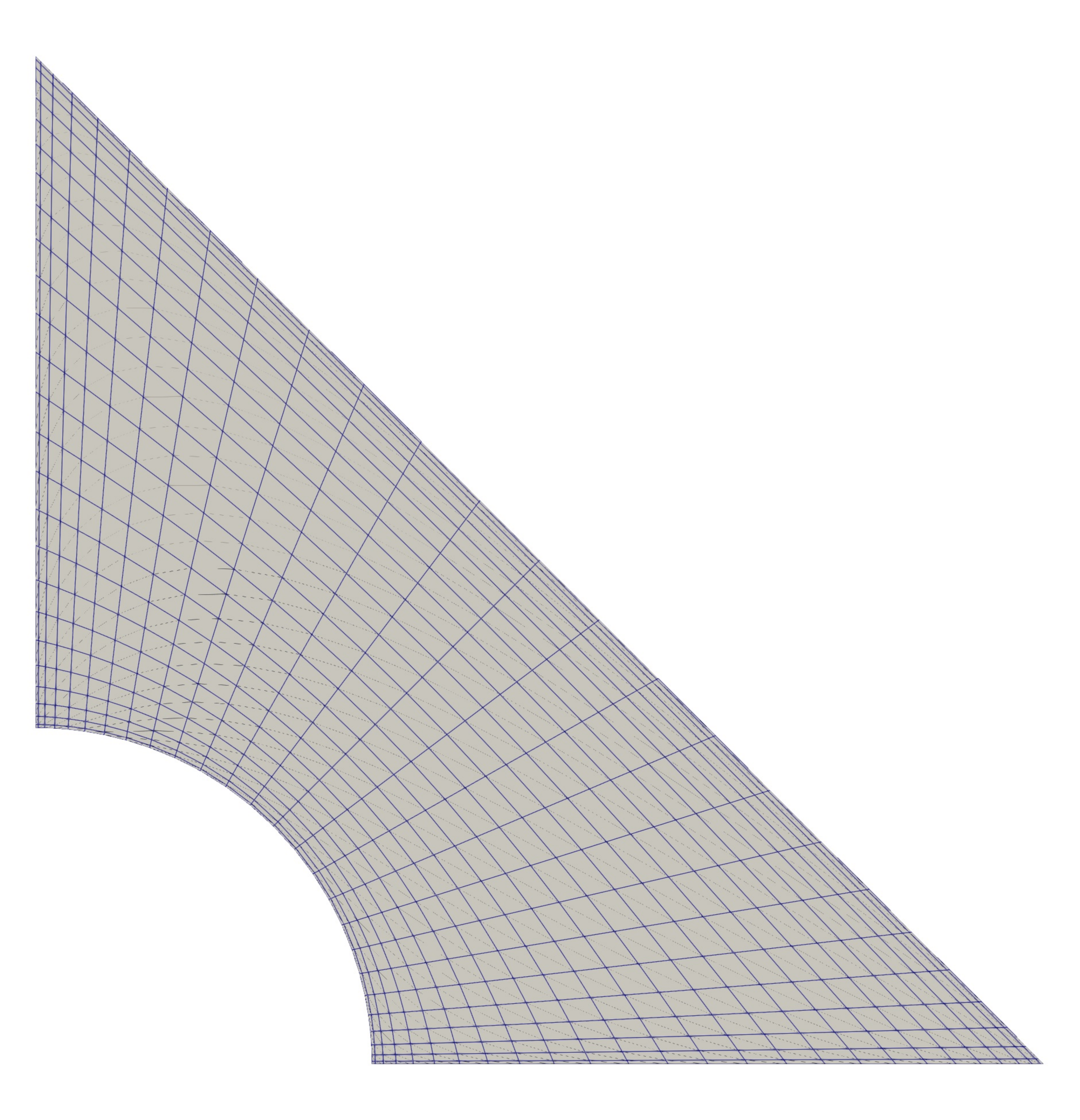} 
    \includegraphics[width=0.495\textwidth]{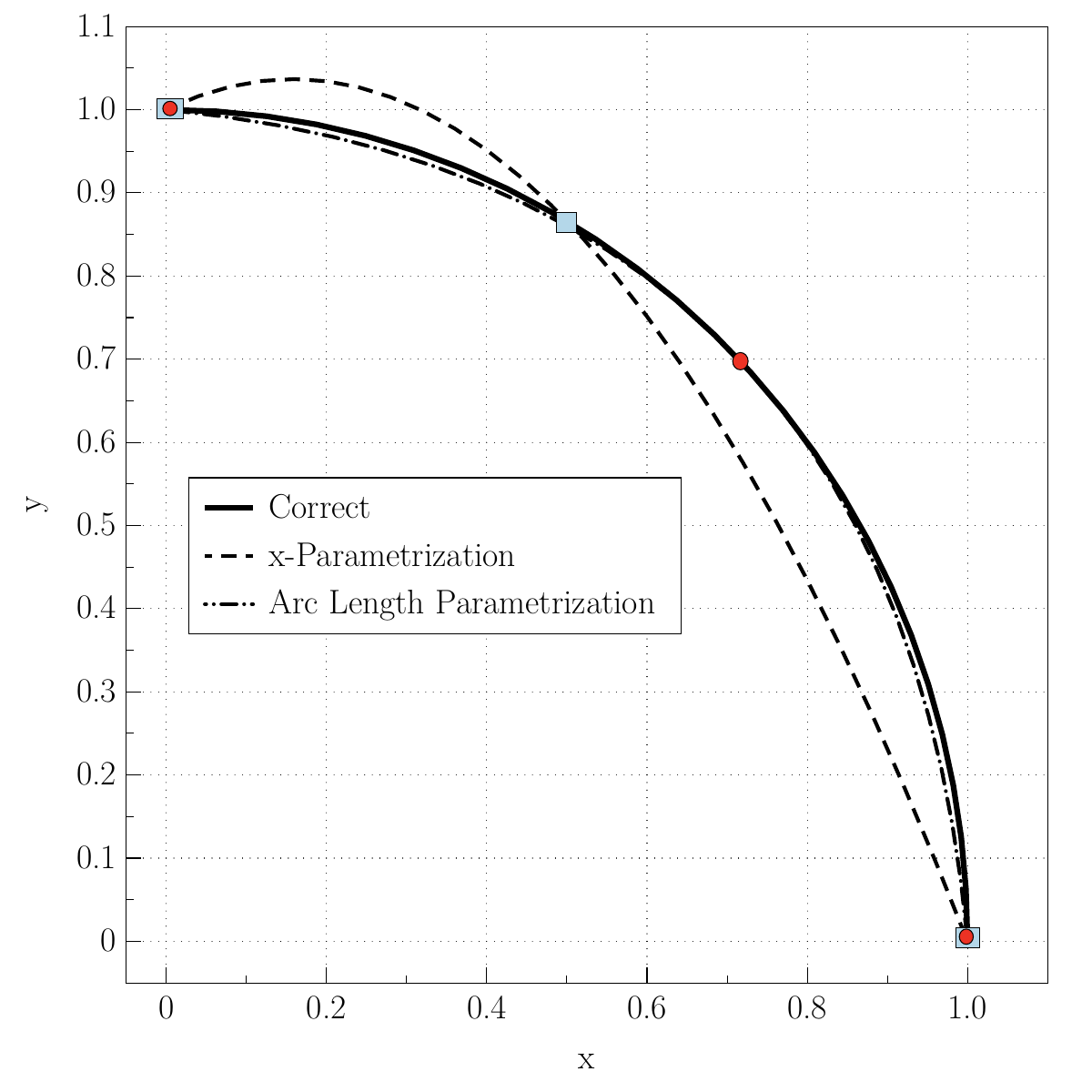} 
   \caption{Left: Correct domain with a circular boundary. Right: Quadratic approximations of the quarter circle using the x-parametrization \eqref{eq:XParametrization} and arc length parametrization \eqref{eq:ArcLengthParametrization}. Squares mark the nodes of the x-parametrization, circles the arc length parametrization}
   \label{fig:CircDomain}
\end{figure}

\begin{figure}[htbp] 
   \centering
   \includegraphics[width=0.495\textwidth]{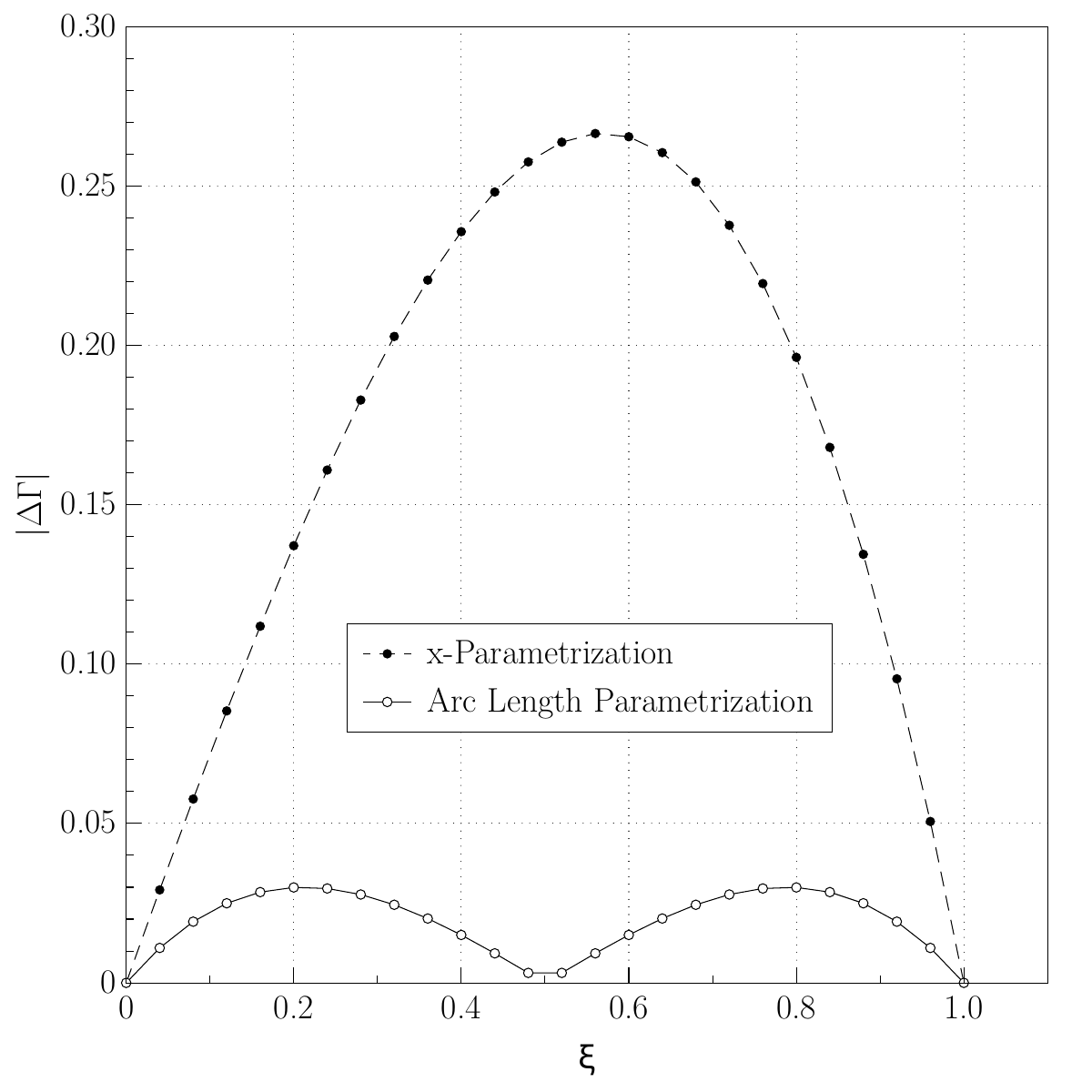} 
   \includegraphics[width=0.495\textwidth]{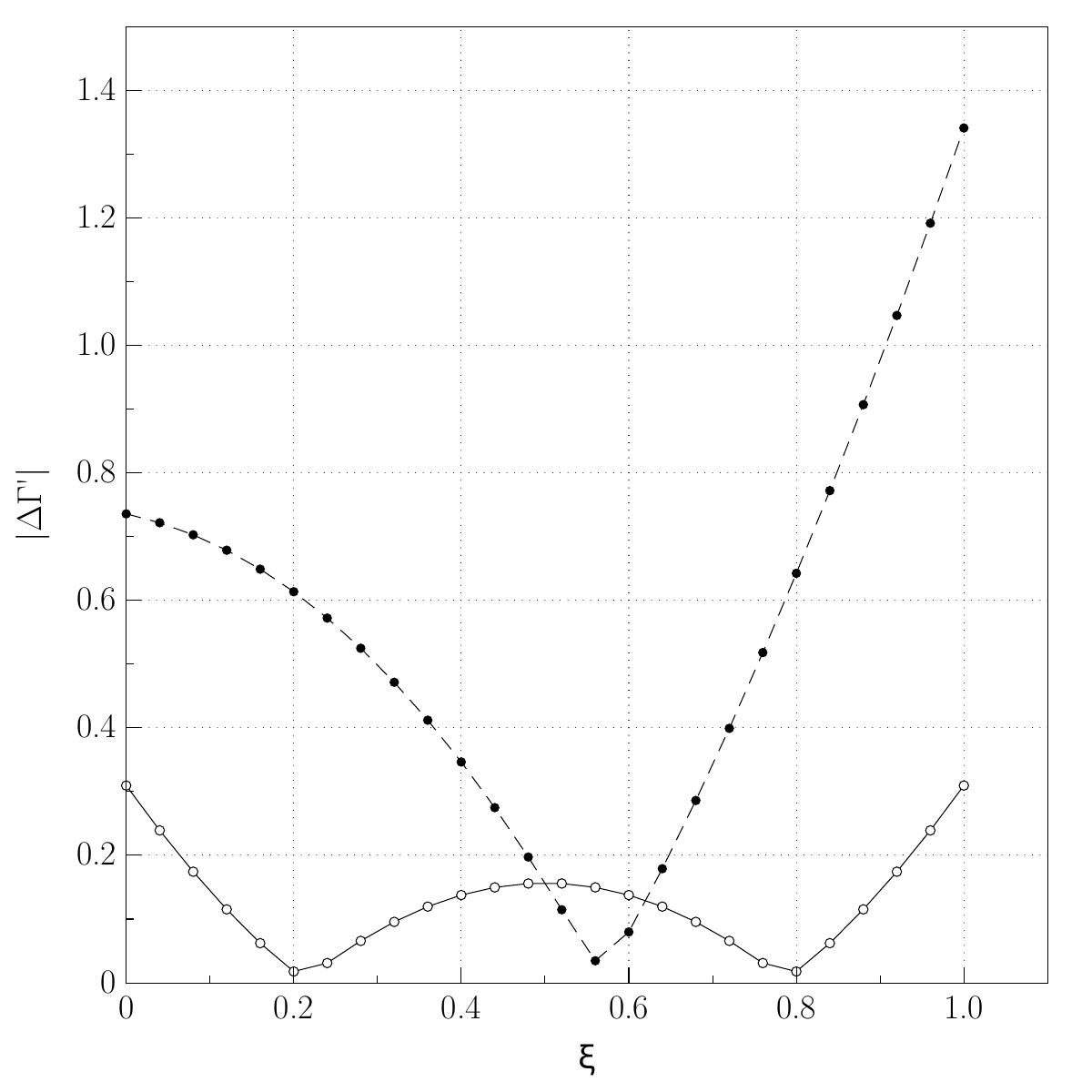} 
   \caption{Errors in location (left) and derivatives (right) for the quadratic polynomial approximation of the circle}
   \label{fig:ParametrizationErrors}
\end{figure}

For the final example we choose a circular arc domain shown in Fig. \ref{fig:CircDomain}, which has a circular domain boundary that will be approximated by quadratic polynomials. It serves to show that how a boundary is approximated, in this case how it is parametrized, affects the error.
We use two parametrizations of the circular boundary. The first is the x-parametrization,
\begin{equation}
\begin{gathered}
x = \xi \hfill\\
y = \sin\left(\frac{\pi}{2}(1-\xi)\right)\hfill
\end{gathered}
\label{eq:XParametrization}
\end{equation}
with the nodes equally spaced placed at $\boldsymbol{\xi} = [0,0.5,1]$. The second is the arc length parametrization,
\begin{equation}
\begin{gathered}
x = \cos\left(\frac{\pi}{2}(1-\xi)\right) \hfill\\
y = \sin\left(\frac{\pi}{2}(1-\xi)\right),\hfill
\end{gathered}
\label{eq:ArcLengthParametrization}
\end{equation}
again, with equally spaced nodes. The differences in the curves are shown in Fig. \ref{fig:CircDomain}.
The errors in location and derivatives are shown in Fig. \ref{fig:ParametrizationErrors}. Table \ref{Tab:MaxErrors} shows the maximum errors and their ratios.

\begin{table}[htp]
\begin{center}
 \caption{Errors for the two parametrizations for the circle\label{Tab:MaxErrors}}
\begin{tabular}{cccc}Approximation & $|\Delta\spacevec \Gamma|_{max}$ & $|\Delta\spacevec \Gamma'|_{max}$ & $||e||_J$ \\
\hline
\hline
x-Parametrization & 0.27 & 1.34 & 0.598 \\
Arc Length & 0.030 & 0.31 & 0.0536 \\
Ratio: & 8.9 & 4.33 & 11.1
\end{tabular}
\end{center}
\end{table}
We solve the same PDE problem as in the previous examples with $N=26$ to get errors less than $1.1\times 10^{-7}$ for both the correct and erroneous domains. The commonly used arc length parametrization gives the better error by an order of magnitude, as seen in Table \ref{Tab:MaxErrors}, which also shows the solution errors and their ratios.

\section{Summary and Conclusions}

Except in the simplest of geometries, the boundaries of meshes on which numerical solutions of PDEs are computed are only {\it approximations} to the true boundaries. The boundary approximation error is one component of the total error in a numerical computation. It is one that is not controlled by the PDE solver, and is especially important for high-order solvers where the solver error could be less than the geometry error. Acknowledgment of this has led to optimization procedures at the mesh generation stage. What has been previously lacking, particularly for hyperbolic systems, is a derivation of the relationship between boundary approximation errors and solution errors.

We have derived the global, $L^2$ error estimates for linear hyperbolic systems generated by errors in the boundary geometry. We show that the error is bounded by data, and is bounded in time if the true solution is bounded in time. 

The total error, \eqref{eq:CharSplitEnergyWGamma}, depends both on volume and boundary errors. The dominant volume errors depend on $\rho = J/J_e = J/(J + \epsilon)$ and $\Delta J\spacevec a^i = J\spacevec a^i_e - J\spacevec a^i$. Along the boundaries, the latter quantities correspond to the errors in the normals. To minimize the volume sources of error in \eqref{eq:CharSplitEnergyWGamma}, then, one would want to minimize the controllable quantities $\epsilon = \Delta J$ and $\Delta J\spacevec a^i$ over the reference domain.
The errors introduced by the boundary conditions appear in the amount of dissipation for outgoing waves, and in the error introduced by evaluating the boundary condition at the wrong location. The boundary location error is second order in the location error, however, whereas the volume error is first order. Therefore, {\it it is not enough to simply evaluate the boundary conditions at the correct locations, even if those are available}, to eliminate the error. Additionally, the results imply that for a numerical computation to converge at the expected rate, the boundary approximations must also converge at least the same order.

We also isolated the effects of a boundary approximation by deriving the error due to a single boundary curve in two space dimensions. There, the error depends to first order on the location error, its derivative, and its second derivative. The numerical experiments here demonstrate this behavior, as do the manufactured solutions in \cite{NORDSTROM2016438}. This suggests that it is desirable to include the error in the normals, through the derivatives of the boundary, as part of an optimization procedure, e.g. as in \cite{ISI:000306588600006,TOULORGE2016361}.

As a final remark, the error as measured by the $L^2$ norm does not address all the important errors in the analysis of hyperbolic systems, where phase errors are also crucially important. The global $L^2$ error, for instance, is not affected by errors at perfectly reflecting (e.g. wall) boundaries, but the normal directions have a critical effect on reflection directions and therefore the locations of nodes and anti-nodes. Thus, there is room for further analysis of the errors due to inaccurate boundary approximation for hyperbolic problems.

\section*{Acknowledgments}

 This work was supported by grants from the Simons Foundation (\#426393,\#961988, David Kopriva).  Andrew Winters was supported by Vetenskapsr{\aa}det, Sweden [award no.~2020-03642 VR]. Jan Nordstr\"{o}m was supported by Vetenskapsr{\aa}det, Sweden [award no.~2021-05484 VR] and the University of Johannesburg Global Excellence and Stature Initiative Funding.

\section*{CRediT authorship contribution statement}

David Kopriva:
Concept, Derivations, Software, Visualization, Writing --- original draft.

Andrew R. Winters:
Concept, Methodology, Writing --- original draft.

Jan Nordstr\"{o}m:
Concept, Methodology, Writing --- original draft.

\section*{Declaration of competing interest}

The authors declare that they have no known competing financial
interests or personal relationships that could have appeared to
influence the work reported in this paper.

\section*{Data availability}

The code used to create the results herein is available upon request
to the corresponding author.

\bibliographystyle{elsarticle-num}
\bibliography{BErrors.bib}

\appendix

\section{Derivation of Error Terms in 2D}\label{AppB}

Here we derive the error term $\mathcal R = -D1 - D2
+\gamma BG
+V1 + V2 + V3 +V4$ in \eqref{eq:CharSplitEnergy2D2} in terms of the geometric errors \eqref{eq:ErrorValues} for the two dimensional domain of Fig. \ref{fig:OneCurved.pdf}.
Those error terms depend on $\nabla_\xi\varepsilon$.  From \eqref{eq:Varepsilon2D},
\begin{equation}
\begin{split}
\nabla_\xi \varepsilon &= \frac{\epsilon \nabla_\xi J - J\nabla_\xi\epsilon}{(J + \epsilon)^2} 
= \frac{\epsilon \nabla_\xi J - J\nabla_\xi\epsilon}{J^2(1 + \epsilon/J)^2}
=  \frac{\epsilon \nabla_\xi J - J\nabla_\xi\epsilon}{J^2}\left(1-2\frac{\epsilon}{J} + O\left(\epsilon^2\right)\right) 
\\&=  \frac{\epsilon\left( \nabla_\xi J +2  \nabla_\xi\epsilon\right)- J  \nabla_\xi\epsilon)}{J^2} + HOTs,
\end{split}
\end{equation}
which needs $\nabla_\xi\epsilon$ and $\nabla_\xi J$. 

First, from \eqref{eq:ErrorValues},
\begin{equation}
\begin{split}
\frac{\partial\epsilon}{\partial\xi} &
= \hat z\cdot\left\{
\left(  \Delta \spacevec \Gamma \times \spacevec X_{\xi\xi}\right)
+\left(  \Delta \spacevec \Gamma' \times \spacevec X_{\xi}\right)
+ (1-\eta)\left[\left(\Delta\spacevec\Gamma''\times \spacevec X_\eta \right) 
+ \left(\Delta\spacevec\Gamma'\times \spacevec X_{\eta\xi} \right) 
- \left( \Delta\spacevec\Gamma''\times \Delta \spacevec \Gamma\right) \right] \right\}
\\&
= \hat z\cdot\left\{
\left( \hat \alpha\times\spacevec X_{\xi\xi} \right)\left|\Delta\spacevec\Gamma\right|
+\left(\hat \beta\times  \spacevec X_{\xi}\right)\left|\Delta\spacevec\Gamma'\right|
+(1-\eta) \left[\left(\hat\gamma\times \spacevec X_\eta \right)\left|\Delta\spacevec\Gamma''\right|
+ \left(\hat \beta\times \spacevec X_{\eta\xi} \right) \left|\Delta\spacevec\Gamma'\right|
- \left( \hat\gamma\times\hat\alpha\right)\left|\Delta\spacevec\Gamma\right|\left|\Delta\spacevec\Gamma''\right|\right]  \right\}
\\&
\equiv {s_1\left|\Delta\spacevec\Gamma\right|+ \left(s_2 + (1-\eta)s_4\right)\left|\Delta\spacevec\Gamma'\right| + (1-\eta)\left[  s_3\left|\Delta\spacevec\Gamma''\right|+ s_5 \left|\Delta\spacevec\Gamma\right|\left|\Delta\spacevec\Gamma''\right|\right]},
\end{split}
\end{equation}
where
\[
s_1 = \hat z\cdot\left( \hat \alpha \times\spacevec X_{\xi\xi}\right),\quad s_2 = \hat z\cdot\left( \hat \beta\times\spacevec X_{\xi} \right)
,\quad s_3 =  \hat z\cdot\left(\hat\gamma\times \spacevec X_\eta \right), \quad s_4 =  \left(\hat \beta\times \spacevec X_{\eta\xi} \right),\quad 
s_5 = \hat z\cdot \left(\hat\alpha\times \hat\gamma\right).
\]
Similarly,
\begin{equation}
\begin{split}
\frac{\partial\epsilon}{\partial\eta} &
 = \hat z\cdot\left\{
\left( \hat\alpha \times \spacevec X_{\xi\eta}\right)\left|\Delta\spacevec\Gamma\right|
- \left(\hat\beta\times \spacevec X_\eta \right) \left|\Delta\spacevec\Gamma'\right|
+ (1-\eta)\left(\hat\beta\times \spacevec X_{\eta\eta} \right)|\Delta\spacevec\Gamma' |
+ \left( \hat\beta\times\hat\alpha\right) \left|\Delta\spacevec\Gamma'\right|\left|\Delta\spacevec\Gamma\right|
 \right\}
\\&
\equiv {p_1\left|\Delta\spacevec\Gamma\right| -\left( p_2 - (1-\eta)p_3\right) | \Delta \spacevec \Gamma'|  +p_4 \left|\Delta\spacevec\Gamma'\right|\left|\Delta\spacevec\Gamma\right|}
 \end{split}
\end{equation}
where
\[
p_1 = \hat z\cdot\left( \hat\alpha\times\spacevec X_{\xi\eta}\right),\quad p_2 = \hat  z\cdot \left(\hat\beta\times \spacevec X_\eta \right)
,\quad p_3 = \hat z\cdot \left(\hat\beta\times \spacevec X_{\eta\eta} \right),\quad p_4 = \hat z\cdot  \left( \hat\beta\times\hat\alpha\right).
\]
Then
\begin{equation}
\begin{split}
J^2 \frac{\partial\varepsilon}{\partial\xi} 
&=
\left( {r_1\left|\Delta\spacevec\Gamma\right|+ (1-\eta)r_2\left|\Delta\spacevec\Gamma'\right| }\right)\left(J_\xi  + 2\left(s_1\left|\Delta\spacevec\Gamma\right|+ \left(s_2 + (1-\eta)s_4\right)\left|\Delta\spacevec\Gamma'\right| + (1-\eta) s_3\left|\Delta\spacevec\Gamma''\right| \right)\right)
\\&
\quad- J\left(  s_1\left|\Delta\spacevec\Gamma\right|+ \left(s_2 + (1-\eta)s_4\right)\left|\Delta\spacevec\Gamma'\right| + (1-\eta)  s_3\left|\Delta\spacevec\Gamma''\right|\right) + HOTs
\\&
={\left(r_1J_\xi - Js_1 \right)\left|\Delta\spacevec\Gamma\right| + \left[ (1-\eta)r_2 J_\xi - J\left(s_2 + (1-\eta)s_4\right)\right]\left|\Delta\spacevec\Gamma'\right| - J(1-\eta)  s_3\left|\Delta\spacevec\Gamma''\right| + HOTs}
\end{split}
\end{equation}
and
\begin{equation}
\begin{split}
J^2 \frac{\partial\varepsilon}{\partial\eta} &= \left( {r_1\left|\Delta\spacevec\Gamma\right|+ (1-\eta)r_2\left|\Delta\spacevec\Gamma'\right| }\right)\left( J_\eta + 2\left(p_1\left|\Delta\spacevec\Gamma\right| -\left( p_2 - (1-\eta)p_3\right) | \Delta \spacevec \Gamma'| \right) \right)
\\&
-J\left(p_1\left|\Delta\spacevec\Gamma\right| -\left( p_2 - (1-\eta)p_3\right) | \Delta \spacevec \Gamma'|   \right) + HOTs
\\&
={\left(r_1 J_\eta -J p_1  \right) \left|\Delta\spacevec\Gamma\right| + \left( (1-\eta)r_2 J_\eta  + J\left( p_2 - (1-\eta)p_3\right)\right)  | \Delta \spacevec \Gamma'| + HOTs}
\end{split}
\end{equation}
where $r_1, r_2,r_3$ are defined in \eqref{eq:riDefs}.

Again from \eqref{eq:ErrorValues}, $\contravec {\mmatrix E}^1\propto \left|\Delta\spacevec \Gamma\right|$ and $\contravec {\mmatrix E}^2\propto \left|\Delta\spacevec \Gamma'\right|$. Therefore, as noted in the observations at the end of Sec. \ref{Sec:ErrorBounds}, $\contraspacevec{\mmatrix E}\cdot\nabla_\xi\varepsilon$, the functional that contributes to $c_1$, has only higher order terms.

On the other hand, the contributor to $c_2$ is
 \begin{equation}
\begin{split}
\contraspacevec{\mmatrix A}\cdot\nabla_\xi\varepsilon &= 
\frac{\contravec {\mmatrix A}^1}{J^2}
\left\{\left(r_1J_\xi - Js_1 \right)\left|\Delta\spacevec\Gamma\right| 
+ \left[ (1-\eta)r_2 J_\xi - J\left(s_2 + (1-\eta)s_4\right)\right]\left|\Delta\spacevec\Gamma'\right| 
- J(1-\eta)  s_3\left|\Delta\spacevec\Gamma''\right|  \right\}
\\&
+ \frac{\contravec {\mmatrix A}^2}{J^2}\left\{ \left(r_1 J_\eta -J p_1  \right) \left|\Delta\spacevec\Gamma\right| + \left( (1-\eta)r_2 J_\eta  + J\left( p_2 - (1-\eta)p_3\right)\right)  | \Delta \spacevec \Gamma'|  \right\} + HOTs
\\&
= { \left|\Delta\spacevec\Gamma\right|\left\{\frac{\contravec {\mmatrix A}^1}{J^2}\left(r_1J_\xi - Js_1 \right) + \frac{\contravec {\mmatrix A}^2}{J^2}\left(r_1 J_\eta -J p_1  \right) \right\} }
\\&
\;{+ | \Delta \spacevec \Gamma'| \left\{\frac{\contravec {\mmatrix A}^1}{J^2} \left[ (1-\eta)r_2 J_\xi - J\left(s_2 + (1-\eta)s_4\right)\right] + \frac{\contravec {\mmatrix A}^2}{J^2} \left[ (1-\eta)r_2 J_\eta  + J\left( p_2 - (1-\eta)p_3\right)\right] \right\}}
\\&
\;{- \left|\Delta\spacevec\Gamma''\right| \left\{\frac{\contravec {\mmatrix A}^1}{J} (1-\eta)  s_3 \right\} + HOTs}
\\&
\equiv{ \left|\Delta\spacevec\Gamma\right|\mmatrix M_1 + | \Delta \spacevec \Gamma'|\mmatrix M_2 + | \Delta \spacevec \Gamma''|\mmatrix M_3 + HOTs},
\end{split}
\end{equation}
where
\begin{equation}
\begin{gathered}
\mmatrix M_1 = \frac{\contravec {\mmatrix A}^1}{J^2}\left(r_1J_\xi - Js_1 \right) + \frac{\contravec {\mmatrix A}^2}{J^2}\left(r_1 J_\eta -J p_1  \right)\hfill\\
\mmatrix M_2 = \frac{\contravec {\mmatrix A}^1}{J^2} \left[ (1-\eta)r_2 J_\xi - J\left(s_2 + (1-\eta)s_4\right)\right] + \frac{\contravec {\mmatrix A}^2}{J^2} \left[ (1-\eta)r_2 J_\eta  + J\left( p_2 - (1-\eta)p_3\right)\right]\hfill\\
\mmatrix M_3 = -\frac{\contravec {\mmatrix A}^1}{J} (1-\eta)  s_3.\hfill
\end{gathered}
\label{eq:Mmatrices}
\end{equation}

Next,
\begin{equation}
\contraspacevec{\mmatrix E}\cdot\nabla_\xi\statevec q = {{{\left|\Delta\spacevec\Gamma\right|}\hat{\mmatrix A}_1}\statevec q_\xi + {\left|\Delta\spacevec\Gamma'\right|(1-\eta)\hat{\mmatrix A}_2}\statevec q_\eta}
\end{equation}
and
\begin{equation}
\varepsilon = {-\frac{1}{J}\left\{r_1\left|\Delta\spacevec\Gamma\right|+ (1-\eta)r_2\left|\Delta\spacevec\Gamma'\right|\right\} + HOTs}.
\label{eq:VarEpsilonExp}
\end{equation}
Then from \eqref{eq:BigCharDefinitions},
\begin{equation}
\begin{split}
V2 &\le
\oneHalf\left|\iprod{ \left|\Delta\spacevec\Gamma\right|\statevec e,\mmatrix M_1\statevec e} \right|
+\left|\iprod{ | \Delta \spacevec \Gamma'|\statevec e,\mmatrix M_2\statevec e}\right|
 + \left|\iprod{ | \Delta \spacevec \Gamma''|\statevec e,\mmatrix M_3\statevec e}\right| + HOTs
\\&
\le 
{ \left\{ \left|\Delta\spacevec\Gamma\right|_{max}\inorm{\frac{\mmatrix M_1}{J}}_{2,\infty}
+   | \Delta \spacevec \Gamma'|_{max}\inorm{\frac{\mmatrix M_2}{J}}_{2,\infty }
+  | \Delta \spacevec \Gamma''|_{max}\inorm{\frac{\mmatrix M_3}{J}}_{2,\infty}\right\}\inorm{\statevec e}_J^2+ HOTs}.
\end{split}
\end{equation}
Next, since $\varepsilon \contraspacevec{\mmatrix E}$ is of higher order,
\begin{equation}
\begin{split}
V3
&\le \left|\iprod{{{\left|\Delta\spacevec\Gamma\right|}\hat{\mmatrix A}_1}\statevec q_\xi, \statevec e} \right|
+\left| \iprod{{\left|\Delta\spacevec\Gamma'\right|(1-\eta)\hat{\mmatrix A}_2}\statevec q_\eta,\statevec e}\right|
\\&
\le {\left|\Delta\spacevec\Gamma\right|_{max}\inorm{\frac{\hat{\mmatrix A}_1\statevec q_\xi}{J^2}}_J\inorm{ \statevec e}_J 
+ \left|\Delta\spacevec\Gamma'\right|_{max}\inorm{\frac{(1-\eta)\hat{\mmatrix A}_2\statevec q_\eta}{J^2}}_{J}\inorm{ \statevec e}_J + HOTs}.
\end{split}
\label{eq:V3AppBBound}
\end{equation}
Finally,
\begin{equation}
\begin{split}
V4 
 &= \iprod{{\frac{1}{J}\left\{r_1\left|\Delta\spacevec\Gamma\right|+ (1-\eta)r_2\left|\Delta\spacevec\Gamma'\right|\right\}\left(\contraspacevec{\mmatrix A}\cdot\nabla_\xi\statevec q\right)} ,\statevec e} + HOTs
 \\&
= \iprod{{\frac{1}{J^2}r_1\left|\Delta\spacevec\Gamma\right|\left(\contraspacevec{\mmatrix A}\cdot\nabla_\xi\statevec q\right)} ,J\statevec e} 
+ \iprod{{\frac{1}{J^2} (1-\eta)r_2\left|\Delta\spacevec\Gamma'\right|\left(\contraspacevec{\mmatrix A}\cdot\nabla_\xi\statevec q\right)} ,J\statevec e} + HOTs
\\&
\le {\left|\Delta\spacevec\Gamma\right|_{max}\inorm{\frac{r_1\left(\contraspacevec{\mmatrix A}\cdot\nabla_\xi\statevec q\right)}{J^3}}_{J}\inorm{ \statevec e}_J 
+ \left|\Delta\spacevec\Gamma'\right|_{max}\inorm{\frac{(1-\eta)r_2\left(\contraspacevec{\mmatrix A}\cdot\nabla_\xi\statevec q\right)}{J^3}}_{J}\inorm{ \statevec e}_J + HOTs}.
 \end{split}
\label{eq:V4AppBBound}
\end{equation}
The quantity $V1$ in \eqref{eq:CharSplitEnergy2D2} is of higher order in the $\Delta$ 's.

What remains is to find the dependencies of the boundary error terms in \eqref{eq:CharSplitEnergy2D2}. The first two boundary terms
represent dissipation error from waves exiting the domain.
Along $\vec \Gamma_1^e$, $\eta = 0$, so with \eqref{eq:VarEpsilonExp},
\begin{equation}
D1 = \int_0^{1}  \frac{1}{J}{ \left\{r_1\left|\Delta\spacevec\Gamma\right|
+ r_2\left|\Delta\spacevec\Gamma'\right| 
\right\} }\statevec e^T \mmatrix A^+\statevec e d\xi + HOTs.
\label{eq:D1}
\end{equation}
Next, from \eqref{eq:ErrorValues},
\begin{equation}
\mmatrix E^{\pm} =\oneHalf \left|\Delta\spacevec\Gamma'\right| \left\{\hat {\mmatrix A}_2 \pm \left| \hat{\mmatrix A}_2 \right|\right\} =  \left|\Delta\spacevec\Gamma'\right| \hat{\mmatrix A}_2^\pm.
\end{equation} 
The term $\varepsilon\mmatrix E^{\pm}$ is second order or higher,
so,
\begin{equation}
D2 
= { \int_{0}^1 \left|\Delta\spacevec\Gamma'\right|\statevec e^T  \hat{\mmatrix  A}_2^+\statevec e d\xi + HOTs.}
\label{eq:BoundarydissipationFnofDelta}
\end{equation}
So to lowest order, the total outflow dissipation error is $O\left(\left|\Delta\spacevec\Gamma\right|,\left|\Delta\spacevec\Gamma'\right|\right)$. 

The next boundary term in \eqref{eq:CharSplitEnergy2D2} is secondary, for it has quadratic dependence on $\left|\Delta\spacevec\Gamma\right|$,
\begin{equation}
BG 
= {\int_{0}^{1} \left|\Delta\spacevec\Gamma\right|^2\left\{\left.(1+\varepsilon)\hat\alpha^T\cdot\nabla_x\statevec g^T \left\{ |\mmatrix A^- | + |\mmatrix E^-|\right\} \nabla_x\statevec g\cdot\hat\alpha\right|_{\eta=0}\right\} d\xi}.
\label{eq:DG}
\end{equation}

When we account for the dissipation being positive and gather terms,
\begin{equation}
\begin{split}
|\mathcal R(\Delta\spacevec\Gamma, \Delta\spacevec\Gamma',\Delta\spacevec\Gamma'')| 
&\le  \left|\Delta\spacevec\Gamma\right|_{max}\left\{\oneHalf \inorm{\frac{\mmatrix M_1}{J}}_{2,\infty } \inorm{ \statevec e}_J 
+\inorm{\frac{\hat{\mmatrix A}_1\statevec q_\xi}{J^2}}_J 
+\inorm{\frac{r_1\left(\contraspacevec{\mmatrix A}\cdot\nabla_\xi\statevec q\right)}{J^3}}_J
\right\}\inorm{ \statevec e}_J
 \\&+ | \Delta \spacevec \Gamma'|_{max}\left\{\oneHalf\inorm{\frac{\mmatrix M_2}{J}}_{2,\infty }\inorm{ \statevec e}_J  
  +\inorm{\frac{(1-\eta)\hat{\mmatrix A}_2\statevec q_\eta}{J^2}}_J +
\inorm{\frac{(1-\eta)r_2\left(\contraspacevec{\mmatrix A}\cdot\nabla_\xi\statevec q\right)}{J^3}}_J
   \right\} \inorm{ \statevec e}_J \\&
 + | \Delta \spacevec \Gamma''|_{max}\left\{\oneHalf\inorm{\frac{\mmatrix M_3}{J}}_{2,\infty }\inorm{ \statevec e}_J  
   \right\} \inorm{ \statevec e}_J + HOTs,
\end{split}
\label{eq:BigRBound}
\end{equation}
so $\mathcal R = O\left( \left|\Delta\spacevec\Gamma\right|, \left|\Delta\spacevec\Gamma'\right|,\left|\Delta\spacevec\Gamma''\right|\right)$.

\end{document}